\documentclass[11pt,a4paper,reqno]{amsart}            

\numberwithin{equation}{section}

\usepackage{amsmath}
\usepackage{amssymb}
\usepackage{graphicx}
\usepackage{color,xcolor}
\usepackage{latexsym}
\usepackage{array, tabularx}
\usepackage{stmaryrd} 
\usepackage{comment}
\usepackage{epsfig}

\usepackage[T1]{fontenc} % pour certains  (ex polonais)
\usepackage[utf8]{inputenc}

\usepackage{soul}

\usepackage[normalem]{ulem}

\def\emm#1,{{\em #1}}

\newcommand{\zs}{\mathbb{Z}}

\newcommand{\qs}{\mathbb{Q}}

\newcommand{\bv}{\bar v}

\newcommand{\bu}{\bar u}

\newcommand{\quasiI}{I_{\scriptscriptstyle{-\!\!\bullet\!\!-}}}

\newcommand{\bm}[1]{\mbox{\boldmath \ensuremath{#1}}}

\newcommand{\beq}{\begin{equation}}
\newcommand{\eeq}{\end{equation}}
\newcommand{\gf}{generating function}
\newcommand{\gfs}{generating functions}

\newcommand{\LandauO}{\mathcal{O}}

\newtheorem{Theorem}{Theorem}[section]
\newtheorem{propo}[Theorem]{Proposition}
\newtheorem{coro}[Theorem]{Corollary}
\newtheorem{Lemma}[Theorem]{Lemma}

\DeclareMathOperator{\Pol}{Pol}

\colorlet{darkgreen}{green!50!black}

\usepackage[plainpages=false,pdfpagelabels,colorlinks=true,citecolor=blue,hypertexnames=false]{hyperref}

\catcode`\@=11
\def\section{\@startsection{section}{1}%
  \z@{.7\linespacing\@plus\linespacing}{.5\linespacing}%
  {\normalfont\bfseries\scshape\centering}}

\def\subsection{\@startsection{subsection}{2}%
  \z@{.5\linespacing\@plus\linespacing}{.5\linespacing}%
  {\normalfont\bfseries\scshape}}

\def\subsubsection{\@startsection{subsubsection}{3}%
  \z@{.5\linespacing\@plus\linespacing}{-.5em}%{.5\linespacing}%
  {\normalfont\bfseries\itshape}}
\catcode`\@=12

\begin{document}

\title[The  degree distribution in bipartite planar maps]
{The  degree distribution in bipartite planar maps:
  applications to the Ising model}

\author[M. Bousquet-M\'elou]{Mireille Bousquet-M\'elou}
\address{M. Bousquet-M\'elou: CNRS, LaBRI, Universit\'e de Bordeaux, 
  351 cours de la Lib\'eration, 33405 Talence, France}
\email{mireille.bousquet@labri.fr}

\author[G. Schaeffer]{Gilles Schaeffer}
\address{G. Schaeffer: CNRS, Laboratoire d'Informatique (LIX), \'Ecole polytechnique, Institut Polytechnique de Paris, 
91128 Palaiseau Cedex, France }
\email{gilles.schaeffer@lix.polytechnique.fr}

\thanks{Both authors were partially
supported by the European Community IHRP
Program, within the Research Training Network "Algebraic Combinatorics
in Europe", grant HPRN-CT-2001-00272.}

%; Secondary}


\begin{abstract}
We characterize the \gf\ of bipartite planar maps counted according to
the degree distribution of their black and white vertices. This result
is applied to the solution of the hard particle and Ising models on
random planar lattices. We thus recover and extend some results
previously obtained by means of matrix integrals.

Proofs are purely combinatorial and rely on the idea that planar maps
are conjugacy classes of trees. In particular, these trees explain  why
the 
solutions of the Ising and hard particle models on maps of bounded
degree are always algebraic.
\end{abstract}

\maketitle

%%%%%%%%%%%%%%%%%%%%%%%%%%%%%%%%%%%%%%%%%%%%%%%%%%%%%%%%
\section*{Foreword}
%%%%%%%%%%%%%%%%%%%%%%%%%%%%%%%%%%%%%%%%%%%%%%%%%%%%%%%%
This paper was first put on arXiv on 4th November 2002~\cite{BMS-arxiv}. Its submission to a journal was then long delayed --- until now. An extended abstract of it was submitted and accepted as a talk at the FPSAC'03 conference (Formal Power Series and Algebraic Combinatorics), held in Vadstena, Sweden, and published in the proceedings of the conference~\cite{BMS-proceedings}.

The reason why we delayed submission is the following: we wanted to find a more combinatorial argument for an enumerative result (Proposition~\ref{propo:Ising4V}, and possibly for generalizations) which is proved using the computation of an integral. 
% establishing [...] in the spirit of .
However, each of us then moved to other projects, and there is no sign that this alternative argument  will ever be written (by us).

Meanwhile, the interest in  this paper has not faded. It has been, and still is,  well cited, either in its arXiv or in its proceedings form; see for instance~\cite{DuSh,Schramm-overview,WW-StFlour,bouttier-mobiles,BDG-blocked,Guionnet-StFlour,AlBo-slice-hypermaps,Sheffield-ICM}. Moreover, the key result of the paper --- a system of algebraic equations describing the \gf\ of planar bipartite maps counted by vertex degrees, see~\eqref{equ:whitetrees}--\eqref{equ:blacktrees}, obtained here via bijective arguments --- has been re-discovered or re-proved by various techniques and various authors, both in  combinatorics and theoretical physics. Let us first cite a paper by Eynard in the physics literature, based on matrix integrals; see~\cite[Sec.~4.3]{Eynard-2matrix} and Eynard's book~\cite[Ch.~8]{eynard-book}. It actually appeared on arXiv, unknown to us, a month before our own preprint, in the High Energy Physics section. Then came an alternative bijective approach by Bouttier, Di Francesco and Guitter, based on mobiles~\cite[p.~21]{bouttier-mobiles}. Let us also cite a very general result by Bernardi and Fusy that includes constraints of the \emm girth, of maps~\cite[Sec.~5.1]{Bernardi-Fusy-hyper}, a  recent new proof based on the fruitful notion of \emm slices,, by Albenque and Bouttier~\cite[Sec.~1.3]{AlBo-slice-hypermaps}. Finally,  an even more recent paper by Albenque, M\'enard and Tokka revisits the bijection of the present paper in terms of orientations~\cite{AMT-Ising}.

Let us return to the reason why we delayed publication, and let us be a bit more technical for a few lines. In this paper, we give combinatorial arguments that explain why the \gf\ of various classes of bipartite maps, counted by their vertex degrees, are algebraic: they are in bijection with certain families of trees, for which algebraicity is somewhat easier to prove.  Our explanations apply mostly to maps in which the root vertex has degree $1$ or degree~$2$. When we are interested in higher root degrees, we resort to an integration procedure, as in Section~\ref{sec:Ising-regular}. From a combinatorics point of view, this is not entirely satisfactory: the antiderivative of an algebraic series need not be algebraic, so that our algebraicity results here look accidental. It is worth mentioning that in several papers where our results are re-proved, such algebraicity results are worked out \emm without any integration procedure,. A central result is of this type is Theorem 8.3.1 in~\cite{eynard-book}, which relies on topological recursion, and is proved by purely combinatorial arguments in~\cite{AlBo-slice-hypermaps}.

Almost 25 years after the completion of our preprint (and a century after the definition of the Ising/Lenz model!), it seems to be high time to submit it to a journal so as to have an officially published version.

%%%%%%%%%%%%%%%%%%%%%%%%%%%%%%%%%%%%%%%%%%%%%%%%%%%%%%%%
\section{Introduction}
%%%%%%%%%%%%%%%%%%%%%%%%%%%%%%%%%%%%%%%%%%%%%%%%%%%%%%%%%
The enumeration of planar maps ({\em i.e.}, connected graphs drawn on 
the sphere with non-intersecting edges) has received a lot of attention
since the sixties. Four decades of exploration have resulted into three
main enumeration techniques, and two types of results: first, a few
remarkable families of maps are counted by very nice and simple
numbers; second, and more generally, many families of maps share the
characteristic of having an {\em algebraic \gf\/}. Meanwhile, many
exciting connections have been established between maps and various
branches of mathematics, like Grothendieck's theory of "dessins
d'enfants" or knot theory, and, most importantly for this paper,
between maps and certain branches of physics.

This paper presents a combinatorial approach for solving two types of
physics models on maps:
 the {\em Ising\/} and {\em hard particle\/} models.  In these models,
the vertices of the maps carry {\em spins\/} or {\em particles\/};
solving these models
is equivalent to a weighted enumeration of 
maps. The weight is in general a specialization of the \emm Tutte
polynomial, of the map~\cite{welsh-book}. 

The idea of putting such a model on maps originates in physics, and
more precisely in two-dimensional {\em quantum gravity\/}: in this
theory, maps,  or
{\em random planar lattices\/}, arise as
discrete models of random geometries.
Br\'ezin \emph{et
al.\/}, in the steps of t'Hooft~\cite{thooft}, developed a new
approach for counting maps, completely different from what had been
done previously  in combinatorics. This approach, suggested by quantum
field theory, is based on the evaluation of matrix integrals
with well chosen potentials~\cite{BIZ,BIPZ}.  An
introduction to these techniques, intended for mathematicians, is
presented in~\cite{Sacha}, and a far reaching account can be found
in~\cite{DFGZJ}.

This approach is extremely powerful: it allows to produce quickly
expressions for generating functions without requiring much invention
at the combinatorial level. In particular, the Ising model on tri- and
on tetravalent maps (maps with vertices of degree 3 or 4) was solved
via this approach in the late 80's, yielding intriguing algebraic
\gfs~\cite{BK87,Ka86,mehta}; the same kind of technique was applied
much more recently to the hard particle model on a larger variety of
maps\footnote{The hard particle model can be seen as a specialization of
the Ising model in a magnetic field, see below for
details.}~\cite{BDFG02a}.
However, the evaluation of the matrix integrals as presented in most
papers is not satisfying from the mathematical point of view:
justifications of the calculations are often omitted, whereas they
involve non-trivial complex analysis and resummation issues.
Moreover, this approach gives little insight on the combinatorial
structure of maps and does not explain the algebraic nature of the
generating functions thus obtained.

\medskip
The material presented in this paper provides an alternative approach
to the Ising and hard particle models, and has none of the above
mentioned drawbacks. It is essentially of a bijective nature, and only
involves elementary combinatorial arguments. We show that certain
families of plane trees are at the heart of both models, and this
justifies the algebraicity of their solutions\footnote{Indeed, what is
more algebraic than a tree \gf ?}. 
We do not study the singularities of the series we obtain, even though
they are  very significant from the physics point of view:
Many results are in this direction are already available in~\cite{BK87,BDFG02a,Ka86};
see also the recent detailed study of~\cite{tokka-Ising-asympt}.

The central objects of our paper are actually \emm bipartite maps,,
which we enumerate according to the degree distribution of their black
and white vertices.  
We thus extend the results of~\cite{BMS00}
dealing with \emm constellations,, and  also
the results of~\cite{BDFG02c} on the degree
distribution of (unicoloured) planar maps 
(these maps are
% the latter being
in one-to-one
correspondence with bipartite maps in which all black  vertices have
degree two).
Then, we show that the solutions of the Ising and hard particle models
are closely related to the enumeration of bipartite maps.
This connection is already used in the physics papers~\cite{BK87,Ka86}.

Our approach provides a new illustration of the principle according to
which planar maps are, in essence, {\em conjugacy classes of trees\/}.
This idea was introduced by the second author of this paper in order
to explain bijectively certain nice formulae for the number of
maps~\cite{Sch97,gillou}. Then, it led to a common generalization of
formulae of Tutte and Hurwitz~\cite{BMS00}, before being adapted to
maps satisfying 2- and 3-connectivity constraints~\cite{domi,PS02}.
Then in 2002, Bouttier \emph{et al.} cleverly built on the same idea
to give a bijective derivation of Bender and Canfield's solution
of the very general problem of counting maps by the degree
distribution of their vertices~\cite{BC94,BDFG02c}
(see~\cite{BMJ06} for the equivalence between both solutions).

Let us mention that the seminal work of Tutte on the enumeration of
maps was not based on bijections, but on recursive decompositions of
maps, which led to functional equations for their \gfs.  This type of
approach is usually fairly systematic; combined with the so-called
quadratic method and its generalizations, it has provided many
algebraicity results, including recent
ones~\cite{BC94,Brown-Tutte,Wormald,Tutte,Tutte-planar}.  It 
yields short self-contained proofs, relying only on elementary
decompositions and algebraic calculus in the realm of formal power
series. However, as far as we know, there has been very few successful
attempts to apply this method to the enumeration of maps carrying an
additional structure, like an Ising model. Roughly speaking, the
standard decompositions of maps still provide certain functional
equations, but these are much harder to solve that in the previous
cases. This is perfectly illustrated by the {\em tour de force\/}
achieved by Tutte, who spent 10 years solving the equation he had
established for the chromatic polynomial of
triangulations~\cite{Tutte-chromatic}.

\medskip For the sake of completeness, let us mention a few additional
results and techniques in map enumeration.  The first bijective proofs 
are due to Cori~\cite{cori} and later to Arqu\`es~\cite{Arq86}.
Some models involving self-avoiding walks on maps can essentially be
solved without matrix integrals~\cite{DK90}, by using enumeration
results due to Tutte.  
Then, a cluster expansion method has been developed
for a general model of spins with neighbor interactions on
maps~\cite{malyshev}.  However, this provides only partial qualitative
results and does not allow to recover the critical exponents of the
hard particle and Ising models (whereas they can be derived from the
\gfs \ obtained either with our approach or with matrix
integrals). Finally, the enumerative theory of ramified coverings of
the sphere has led a number of authors to consider a very refined
enumeration of maps on higher genus surfaces, in which the degrees of
vertices {\em and\/} faces are taken into account. The relevant
methods then rely on the encoding of maps by permutations~\cite{CM}
and lead to expressions involving characters of the symmetric
group~\cite{GJ97,JV90b}. 

The structure of the paper is outlined at the end of the next section, after the description of our results.
 
%%%%%%%%%%%%%%%%%%%%%%%%%%%%%%%%%%%%%%%%%%%%%%%%%%%%%%%%%%%%%%
\section{A glimpse at the results}
\label{sec:glimpse}
\begin{figure}[htb]
\begin{center}
\epsfig{file=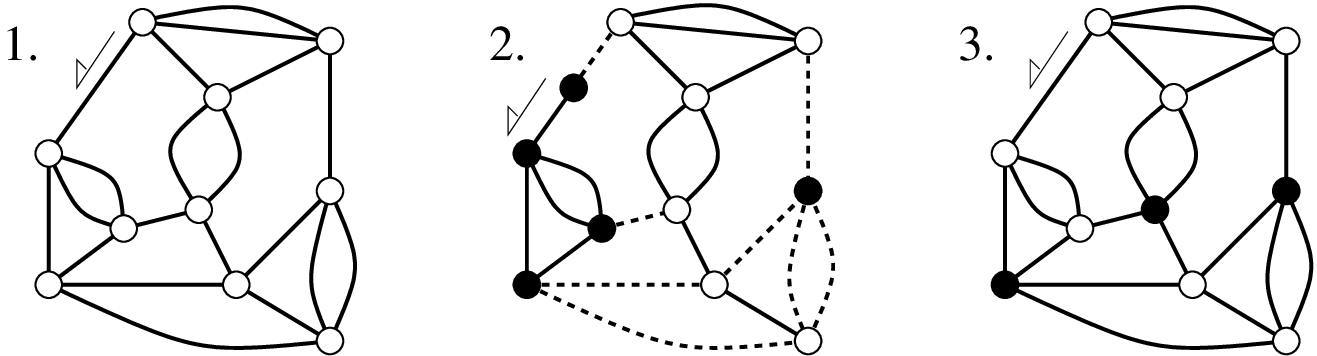, width=12cm}
\end{center}
\caption{1. A tetravalent map -- 2. An Ising configuration on a
  quasi-tetravalent map with 8 frustrated edges (dashed edges) -- 3.  A
  hard particle configuration with 3 particles on a tetravalent map. }
\label{fig:example-ising}
\end{figure}

%%%%%%%%%%%%%%%%%%%%%%%%%%%%%%%%%%%%%%%%%%%%%%%%%%%%%%%%
To begin with, let us recall a few definitions and conventions. A
planar map is a 2-cell decomposition of the oriented sphere into
vertices (0-cells), edges (1-cells) and faces (2-cells). In more
vernacular terms, it is a connected multigraph drawn on a sphere with
non-intersecting edges. Loops and multiple edges are allowed. Two maps
are isomorphic if there exists an orientation preserving homeomorphism
of the sphere that sends one onto the other.
The  degree of a vertex is the number of incident half-edges, which means that a loop counts twice. The  degree of a face is the length of its boundary, where each edge has length $1$.
A map is {\em rooted\/} if one of its edges is distinguished and
oriented.  In this case, the map is drawn on the plane in such a way
the infinite face lies to the right of the root edge.  All the maps
considered in this paper are planar, rooted, and considered up to
isomorphisms.  An example  is provided in
Figure~\ref{fig:example-ising}.1; this map is tetravalent, meaning that
all vertices have degree $4$. In the physics literature, authors
often consider unrooted maps; but they count them with a symmetry
factor which makes the problem equivalent to counting rooted maps, up
to a differentiation of the (unrooted) \gf; see e.g.~\cite[Sec.~5.2.3]{handbook-schaeffer} for details.

%============================================================
\subsection{The Ising and hard particle models on maps}
\label{subsec:Ising-hard}
% ============================================================

In combinatorial terms, the solution of the Ising model on planar maps
(in a magnetic field) is equivalent to counting maps having vertices
of two types (white and black): one has to count them, not only by
their number of edges and vertices of each colour, but also by the
number of {\em frustrated edges\/}, that is, edges that are adjacent
to a white vertex and to a black one.

The tools developed in this paper allow us to solve the Ising model on
any class of maps subject to degree conditions. The associated \gfs\ 
are algebraic as soon as the degrees of the vertices are bounded. Here
is, for instance, the result that we obtain for ``quasi-tetravalent'' maps:
maps having only vertices of degree $4$, except for the root vertex
which has degree $2$
(Figure~\ref{fig:example-ising}.2).

\begin{propo}[Ising on quasi-tetravalent maps]\label{propo:Isingquasi4V}
  Let $\quasiI(X,Y,u)$ be the \gf\ of bicoloured quasi-tetravalent maps,
  rooted at a black vertex, where the variable $X$ (resp.~$Y$) counts
  % the number of 
  white (resp.~black) vertices of degree $4$, and $u$
  counts
  % the number of
  frustrated edges. Let $P(x,y,v)\equiv P$ be the power
  series defined by the following algebraic equation:
  \[
  P=1+3xyP^3 + v^2 \ \frac{P(1+3xP)(1+3yP)}{(1-9xyP^2)^2}.
\]

  Then the Ising \gf \ $\quasiI(X,Y,u)$ can be expressed in terms of
  $P(x,y,v)$, with $x=X(u- \bu)^2, y=Y(u-\bu )^2$, and $v= \bu =
  1/u$. One possible expression is:
\[
\frac{\quasiI(X,Y,u)}{1-\bu^2}=xP^3+\frac{P(1+3xP-2xP^2-6xyP^3)}{1-9xyP^2}
- \frac{y v^2P^3(1+3xP)^3}{(1-9xyP^2)^3}.
\]
As $P$ itself, the series $\quasiI(X,Y,u)$ is algebraic of degree $7$.
\end{propo}

\begin{figure}[htb]
\begin{center}
\epsfig{file=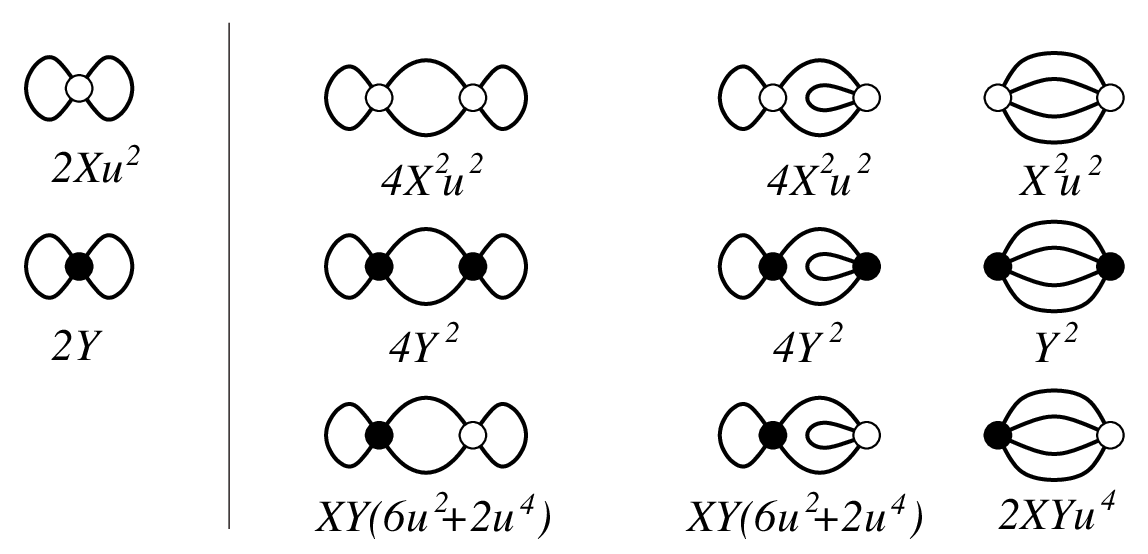, width=12cm}
\end{center}
\caption{The (unrooted) tetravalent maps on the sphere having
   one or two vertices. The polynomials
  %in $X,Y$ and $u$ shown in the figure
 describe
 %the number of ways of adding
 the total Ising weight of
the  quasi-tetravalent maps obtained by adding a black root vertex of degree $2$ on an edge and rooting the resulting map at one of the two incident half-edges.}
\label{fig:FirstOnes}
\end{figure}

\noindent {\bf Remarks} \\
{\bf 1.}
Replacing $X$ by $ t X$  and $Y$ by $ t Y$ 
gives a \gf \ in
which tetravalent vertices are counted by $t$. The expansion of the
Ising \gf \ then  begins as follows: 
\[
\quasiI( tX, tY,u) =1+ t(2Xu^2+2Y)+t^2\left(9X^2u^2+9Y^2+XY(12u^2+6u^4)\right)
+\LandauO(t^3).
\]
These terms correspond to maps having at most two tetravalent vertices,
% they are
shown in Figure~\ref{fig:FirstOnes}. \\
{\bf 2.}
It is easy to check that the
parametrization by $P$ is equivalent to the one given by Boulatov
and Kazakov for the free energy of the Ising model  on tetravalent
maps~\cite[Eq.~(17)]{BK87}. 

\medskip 
Proposition~\ref{propo:Isingquasi4V} 
will only be proved in Section~\ref{sec:ising}. 
At the heart of the proof is the fact that the series  $P$ is the \gf \
of certain  trees, as suggested 
by  the equation defining~$P$.  We shall also prove
 that the Ising \gf\ of truly tetravalent maps belongs to the
same algebraic extension of $\qs(X,Y,u)$ as $P$. But its expression in
terms of $P$ is messier (Proposition~\ref{propo:Ising4V}).

Observe that the series $\quasiI(X,Y,u)$ contains all the information that we
need to count also the number of {\em monochromatic\/} (i.e., non-frustrated)
edges.
Depending on the colour of their endpoints, we call them white edges or black edges. In particular, the series 
$u\quasiI(X,Yu^2,1/\sqrt u)$ 
counts bicoloured quasi-tetravalent maps by their white and black vertices
(variables $X$ and 
$Y$), and by the number of  black edges (variable $u$). By
setting $u$ to zero in this series, we forbid such edges.
In particular, both neighbours of the black root vertex are white. 
Let us erase the root vertex: the root edge is now  white.
Let $H(X,Y) $ denote the limiting series of $u\quasiI(X,Yu^2,1/\sqrt u)$ as $u$ goes
to zero: it counts bicoloured planar maps, rooted at a
 white edge, in which two adjacent vertices cannot be both
black. Say that a white (resp. black) vertex is \emm vacant, (resp. \emm occupied
by a particle,): we have just
solved the so-called {\em hard particle model\/} on tetravalent maps.
A hard particle configuration is
 shown in Figure~\ref{fig:example-ising}.3.

 \begin{coro}[Hard particles on tetravalent maps]
   \label{cor:hard4V}
The hard particle \gf\ of tetravalent maps rooted at an edge whose
endpoints are vacant is algebraic of degree $7$ and can be expressed as:
\[
H(x,y)=xP^3+\frac{xP^2(3-2P)}{1-9xyP^2} -\frac{27x^3yP^6}{(1-9xyP^2)^3}
\]
where $P\equiv P(x,y)$ is the power series defined by
\[
P=1+ 3xyP^3 + \frac {3xP^2}{(1-9xyP^2)^2}.
\]
\end{coro}

More generally, our techniques will allow us to solve the hard
particle model on any class of maps defined by degree conditions. This
includes the case of tri- and tetravalent maps solved in~\cite{BDFG02a} via
matrix integrals, but not the case of tri- and tetravalent {\em
bipartite} maps  (all cycles have an even length) also solved
in~\cite{BDFG02a}. These bipartite models seem to mimic more faithfully 
the phase transitions observed on the square and honeycomb
lattices (which are of course bipartite).

% ============================================================
\subsection{Bipartite maps and blossom trees}
\label{sec:bip-init}
% ============================================================
Our method for solving the Ising and hard particle models uses a
detour via the enumeration of {\em fully frustrated\/} maps, better
known as {\em bipartite\/} maps. Their enumeration will be based on a
correspondence between these maps and certain trees, called {\em
blossom trees\/}. 
These trees have an algebraic \gf \ as soon as the
degrees of their vertices are bounded. In particular, the series $P$
of Proposition~\ref{propo:Isingquasi4V}
% and Corollary~\ref{cor:hard4V}
will be shown to count certain blossom trees.

\begin{figure}[htb]
\begin{center}
\epsfig{file=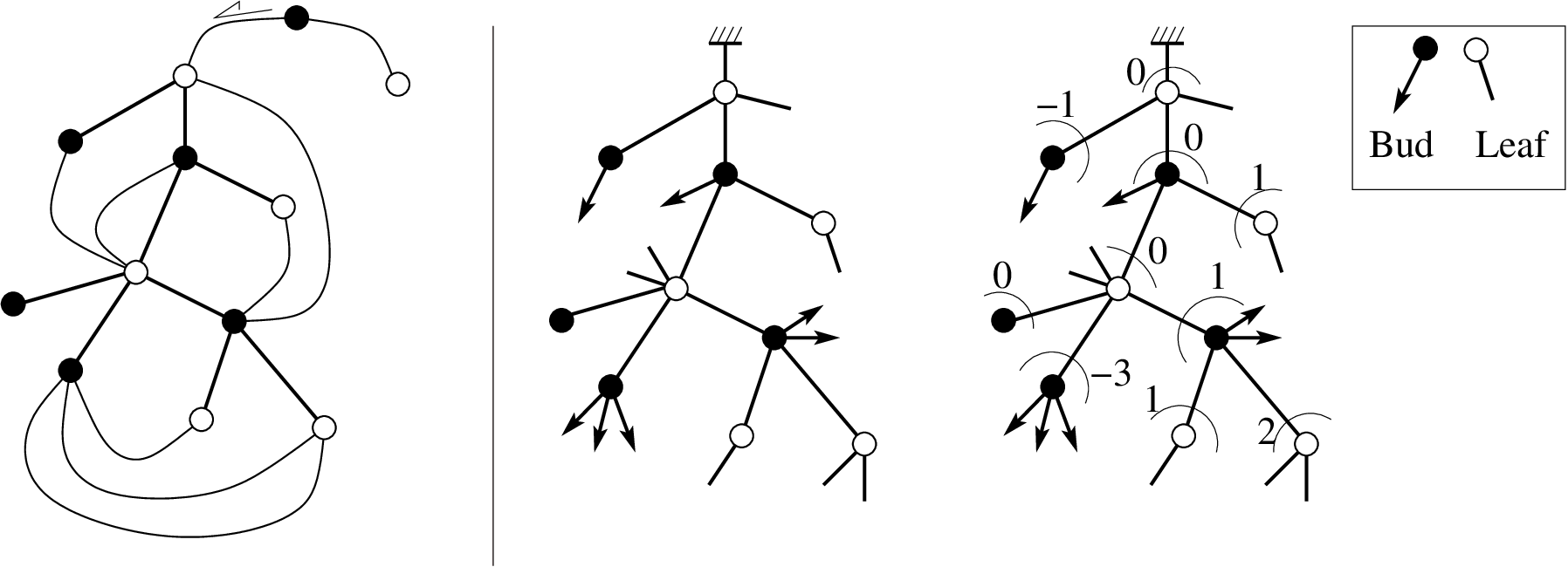, width=14cm}
\end{center}
\caption{A bipartite map; a blossom tree rooted at a leaf, 
and the charge at each subtree.}
\label{fig:example}
\end{figure}

Let us give a brief description of these blossom trees.
% that are at the heart of bipartite maps.
An example is provided by Figure~\ref{fig:example}.
These trees are one-face maps, and are \emm properly,
  % As one can expect, these trees are themselves
  bicoloured: all the neighbours of a black vertex are white,
and vice-versa.  In addition to the standard vertices and edges,
blossom trees carry {\em half-edges\/}, which are called {\em
leaves\/} when they hang from a white vertex, and {\em buds\/}
otherwise. 
Leaves are represented in our figures by short segments, while {buds}
are represented  by black outgoing arrows.
  The trees are rooted at a leaf or a bud, and the vertex
attached to this half-edge is called the root vertex.
Finally, we define the {\em charge\/} of a tree to be the
difference between the number of leaves and the number of buds it
contains; the root half-edge does not count. The charge at a vertex is
the charge of the subtree rooted 
at this vertex. We impose the following charge conditions on the
\emm non-root, vertices of blossom trees:
%(except at the root vertex): 
white vertices have a nonnegative 
charge, while %all
black vertices have charge at most one.  The
notion of charge was introduced by Bouttier {\em et al.\/}
in~\cite{BDFG02c} for the enumeration of usual planar maps, and turns
out to be also useful in the bipartite case.

Given the neat recursive structure of trees, it is not difficult to
write functional equations that govern a very detailed \gf \ of
blossom trees. Let $\bm x = (x_1, x_2, \ldots)$ and $\bm y = (y_1,y_2,
\ldots)$. For $i \in \zs$, let $W_i(\bm x,\bm y)$ be the \gf \ of
blossom trees rooted at a leaf such that the charge at the (white)
root vertex is $i$. In this series, the variable $x_k$ (resp.~$y_k$)
counts white (resp.~black) vertices of degree $k$. Similarly, 
let $B_i(\bm x,\bm y)$ be the \gf\ of blossom trees rooted
at a bud such that the charge at the (black) root vertex is $i$.
 Let us form the following  \gfs:
\beq
\label{WBdef}
W(\bm x,\bm y;z) \equiv W(z)= \sum_{i \ge 0 } W_i(\bm x,\bm y)\, z^i
\quad\hbox{ and }\quad
\quad B(\bm x,\bm y;z) \equiv B(z) = \sum_{i \le 1} B_i(\bm x,\bm y)\, z^i .
\eeq
Then
\begin{equation}\label{equ:whitetrees}
W(z) = [z^{\ge 0}] \sum_{k\ge 0} x_{k+1} \Big(z+ B(z)\Big)^k,
\end{equation}
and
\begin{equation}\label{equ:blacktrees}
B(z) = [z^{\le 1}] \sum_{k\ge 0} y_{k+1} \left( \frac 1 z +
  W(z)\right) ^k.
\end{equation}
We have used the following notation: for any power series $S$ in $\bm
x $ and $\bm y $ having coefficients in $\zs[z, 1/z]$, the series
$[z^{\ge 0}]S(\bm x,\bm y;z)$ is obtained by selecting from $S$ the
terms with a nonnegative power of $z$.  The notation $[z^{\le 1}]$
% naturally
corresponds to the extraction of terms in which the exponent
of~$z$ is at most one. 
More generally, for any  $i \in \zs$,
\begin{equation}\label{equ:blackgf}
  B_i (\bm x,\bm y) = [z^i] \sum_{k \ge 0} y_{k+1} \left( \frac 1 z +  W(z)\right) ^k
  \quad  \text{and}   \quad
  W_i (\bm x,\bm y) = [z^i] \sum_{k \ge 0} x_{k+1} \Big( z +  B(z)\Big) ^k.
\end{equation}

\medskip
Our central theorem gives an expression for the 
 {\em degree \gf \/} of
bipartite planar maps. This series, denoted by $B(\bm x,\bm y)$,
enumerates maps by their number of white and black vertices having a
given degree $k$ (variables $x_k$ and $y_k$). The theorem below thus
solves a combinatorial problem that has an interest of its
own. Moreover, we shall see that the solutions of the Ising and hard
particle models on planar maps have close connections to it.

\begin{Theorem}\label{thm:main}
The degree \gf \ of bipartite planar maps rooted at a black vertex of
degree $2$ is related to the \gfs\ of blossom trees as follows:
\[
  M(\bm x,\bm y) = y_2 \left((W_0-B_2)^2+ W_1 -B_3 -B_2^2 \right),
\]
where $W_i\equiv W_i(\bm x, \bm y)$ and  $B_i\equiv B_i(\bm x, \bm y)$
are defined by~\eqref{WBdef}--\eqref{equ:blackgf}.
\end{Theorem}

 This theorem~\ref{thm:main} will be proved in  Section~\ref{sec:conjugue}.
 
\begin{figure}[tb]
\begin{center}
\epsfig{file=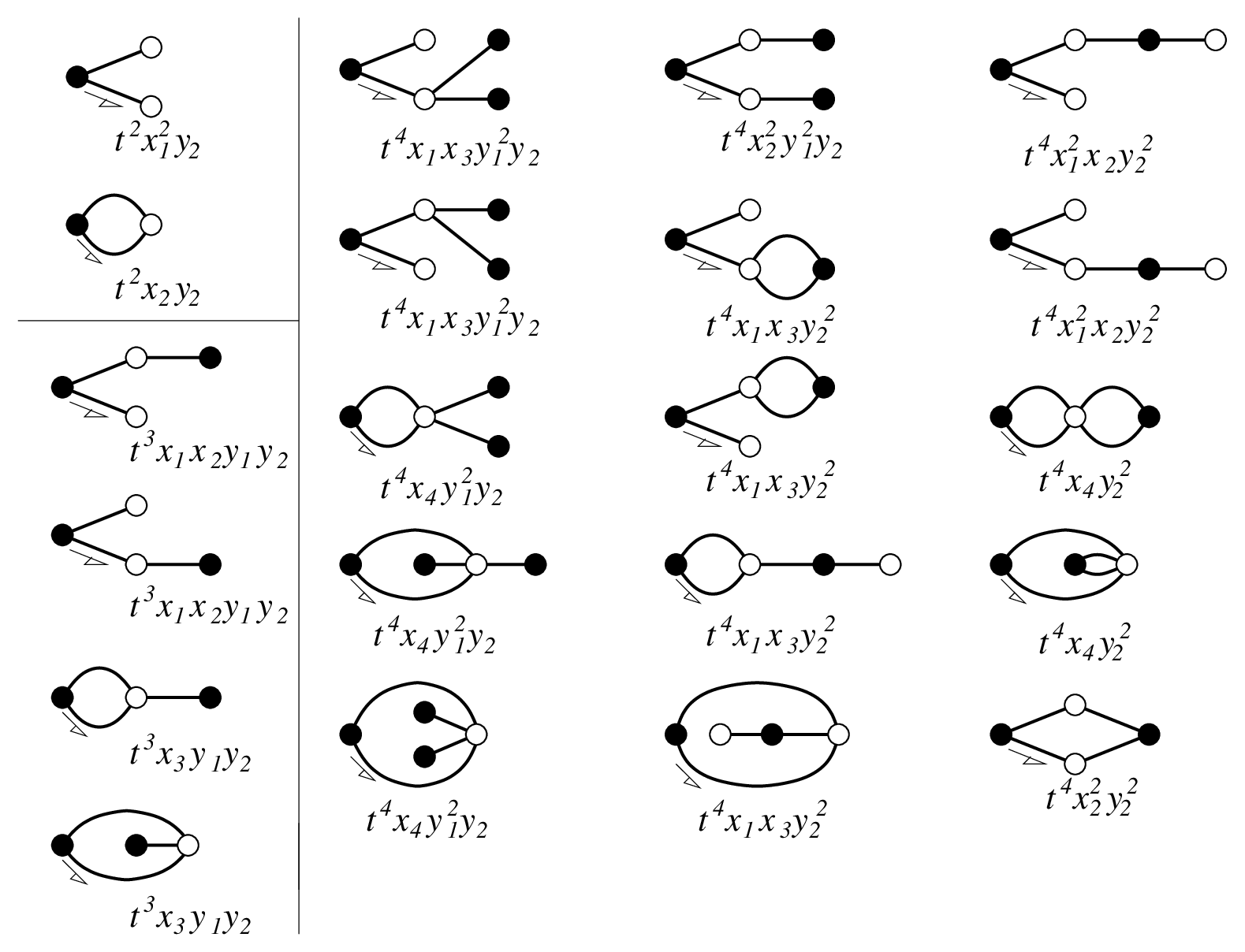, width=14cm}
\end{center}
\caption{The bipartite maps  rooted on a black vertex of degre two
  and having at most $4$ edges, together with their degree
  distribution.}
\label{fig:bipartite}
\end{figure}

\medskip
\noindent{\bf Remarks}\\
{\bf 1.} One can of course write  $M(\bm x,\bm y) = y_2 (W_0^2+ W_1 -B_3
-2W_0B_2 )$, and this expression is better suited to the
applications. However, the expression of Theorem~\ref{thm:main} will be shown
to reflect more faithfully the combinatorics of maps.\\
{\bf 2.} 
Replacing $x_k$ by $t^k x_k$  gives a \gf \ in
which edges are counted by $t$. The first few terms of $M$ then read:
\[
t^2y_2(x_1^2+x_2)+2t^3y_2y_1(x_1x_2+x_3)
+t^4y_1^2y_2(2x_1x_3+3x_4+x_2^2)+t^4y_2^2(4x_1x_3+2x_1^2x_2+2x_4+x_2^2)
+\LandauO(t^5).
\] 
The corresponding bipartite maps are shown in
Figure~\ref{fig:bipartite}.\\
{\bf 3.} Let us rephrase the above result in terms of
factorizations of permutations. A bipartite planar map with~$n$  edges labelled from $1$ to $n$ can be
encoded by a pair of permutations $(\sigma,\tau)$ of the symmetric group $\mathcal{S}_n$,
satisfying the following conditions: the group generated by $\sigma$
and~$\tau$ acts transitively on $\{1,\ldots,n\}$, and the three
permutations $\sigma$, $\tau$ and $\phi:=\sigma\tau$ have a total of $(n+2)$
cycles~\cite{CM}.  The enumeration of such \emph{minimal transitive
factorizations} is the subject of a vast litterature; see~\cite{BMS00,
GJ97} and references therein.
For $\lambda$ and $\mu$ two partitions of~$n$, let $m(\lambda,\mu)$ be
the number of pairs of permutations $(\sigma,\tau)$ of respective
cyclic type $\lambda$ and $\mu$, satisfying the two conditions
indicated above. Then the series of Theorem~\ref{thm:main} can be written as
% gives an expression for the generating function
\[
  M(\bm x,\bm y)
=\sum_{n\geq1}\sum_{\lambda,\mu}m(\lambda,\mu)\cdot 2m_2\cdot
\bm{x}^{\lambda}\bm{y}^{\mu}/n!,
\]
where the inner summation is on all
partitions $\lambda=1^{\ell_1}\ldots n^{\ell_n}$ and
$\mu=1^{m_1}\ldots n^{m_n}$ of $n$, and the weights
$\bm{x}^\lambda\bm{y}^{\mu}$ are defined by
$\bm{x}^\lambda=x_1^{\ell_1}\cdots x_n^{\ell_n}$, and
$\bm{y}^\mu=y_1^{m_1}\cdots y_n^{m_n}$.

\bigskip

\noindent{\bf Outline of the paper.}
% The paper is organized as follows:
In Section~\ref{sec:closure}, we
describe a general connection between maps and trees: one transforms a
map into a tree by cutting %certain
some of its edges into two half-edges.
Conversely, merging half-edges of a tree gives a planar map.  In
Section~\ref{sec:bijection}, we show that these transformations,
restricted to certain classes of maps and trees (called balanced
trees) are one-to-one. Balanced trees are then enumerated in
Section~\ref{sec:conjugue}: this proves our central
Theorem~\ref{thm:main} above. Finally, Sections~\ref{sec:hard}
and~\ref{sec:ising} are respectively devoted to applications of this
theorem to the solution of the hard particle and Ising models.

%%%%%%%%%%%%%%%%%%%%%%%%%%%%%%%%%%%%%%%%%%%%%%%%%%%%%%%%%%%%%%%
\section{Maps and trees}
\label{sec:closure}
%%%%%%%%%%%%%%%%%%%%%%%%%%%%%%%%%%%%%%%%%%%%%%%%%%%%%%%%

%============================================================
\subsection{Maps with legs and blossom trees}
\label{sec:legs}
% ======================================

\begin{figure}[htb]
\begin{center}
\epsfig{file=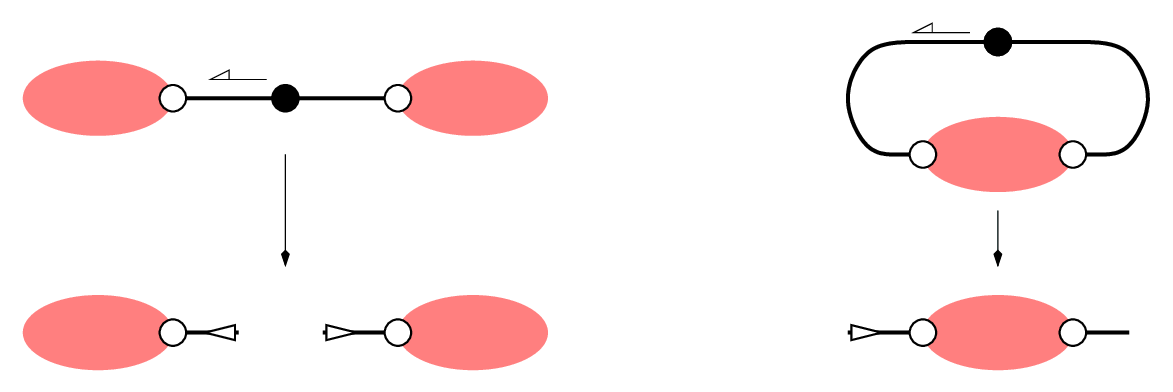, width=13cm}
\end{center}
\caption{The deletion of a black root vertex of   degree $2$.}  
\label{fig:rootedortwo-legs}
\end{figure}

Take a bipartite planar map $M$ rooted at a black vertex $v$ of degree
$2$. Such maps are exactly those we wish to count
(Theorem~\ref{thm:main}). 
Let us consider that each of the two edges
% that start from
incident to~$v$ is made of two half-edges. Delete $v$ and the two half-edges
attached to it. Two cases occur (Figure~\ref{fig:rootedortwo-legs}):
\begin{itemize}
\item either we obtain an ordered pair of maps, each of them carrying a
half-edge (or {\em leg\/}), on which it is rooted,
\item or we obtain a single map with two half-edges, and root this map at
the half-edge belonging to the root edge of $M$.
\end{itemize}

The degree \gf \ of bipartite maps rooted at a black vertex of degree
 $2$ can thus be written as \beq
\label{equ:maps-pattes} 
M(\bm x, \bm y) = y_2( L_1(\bm x, \bm y)^2 + L_2(\bm x, \bm y)), 
\eeq 
where the
series $L_1$ and $L_2$ respectively count by their degree distribution
maps with one or two legs, rooted at a leg. The reason why we keep the
legs is that we do not want to modify the degrees of the vertices.
Compare Equation~\eqref{equ:maps-pattes} to Theorem~\ref{thm:main}:
% the next three sections
this section and the next two
are devoted to proving that $L_1=W_0-B_2$ and $L_2=
W_1-B_3-B_2^2$.

The above decomposition explains why, as in~\cite{BDFG02c}, we shall be
interested in generalized bipartite maps that do not only consist of
the %traditional
usual edges and vertices, but also contain a number of
half-edges in their infinite face (when drawn on the plane), and are
rooted at one of these half-edges.  Let us insist on the distinction
between a half-edge, which is incident to one vertex only, and a
\emm dangling, edge, which is incident to two vertices, one of which has
degree one.  \emph {From now on, in this section and the next two, all
  maps are bipartite (with a proper bicolouring)  and rooted at a half-edge.}

\medskip
We need to introduce more terminology. In passing, we shall
reformulate slightly the definition  of blossom trees given in
Section~\ref{sec:glimpse}.
The vertex adjacent to the root half-edge is called the root vertex of
the map. Half-edges that hang from black vertices are called
\emph{buds} and are represented in our figures by black outgoing
arrows.  Half-edges that hang from white vertices are called
\emph{leaves} and are represented by short segments.  The \emph{degree
  distribution} of a map is the pair of partitions $(\lambda,\mu)$
such that $\lambda$ gives the degree distribution of white vertices
and $\mu$ gives the degree distribution of black vertices. Half-edges
are included in the degree of the vertex they are attached to.  A map
with $n_b$ buds and~$n_\ell$ leaves %obviously
satisfies
$n_b+|\lambda|=n_\ell+|\mu|$ where $|\lambda|$ denotes the sum of the
parts of $\lambda$.  The $1$-leg map of Figure~\ref{fig:closure} has degree distribution $\lambda= 2^2 \cdot 3\cdot4\cdot5$, $\mu=2\cdot 4^2\cdot5$.

Let us now define two important subclasses of maps. A \emph{$k$-leg
map} is a map with exactly $k$ leaves and no bud; hence a $k$-leg map
is  rooted at one of its leaves.  A \emph{tree} is a map with
only one face (and an arbitrary number of buds and leaves),  rooted at a half-edge as all maps in this section. In our figures, the root is usually at the top. The
\emph{total charge} of a tree is the difference between its number of
leaves and its number of buds; the root half-edge is counted. The
\emph{charge} of a tree is the same difference, but the root half-edge
is not counted. Hence the charge and total charge always differ by
$\pm 1$.

Take a tree $T$ and an edge $e$ of this tree. Cut $e$ into two
half-edges: this leaves two {subtrees}, rooted at these half-edges. 
The subtree that does not contain the root of $T$ is called the \emph{lower subtree} of $T$ at $e$. 
Let $T^\bullet_e$ denote the subtree containing the black
endpoint of $e$ and $T^\circ_e$ the other one.
%{We root them at the half-edge coming from $e$.} 
%
The charges $c_\circ$ of $T_e^\circ$ and $c_\bullet$ of $T_e^\bullet$
satisfy $c_\circ+c_\bullet=c$, where $c$ is the total charge of~$T$.
The \emph{black charge rule} is satisfied  at $e$ if the subtree
$T_e^\bullet$ has charge at most one.  The \emph{white charge rule} is
satisfied  at $e$ if the subtree $T_e^\circ$ has a nonnegative charge.

A tree is $T$ called a \emph{blossom tree} if all its \emph{lower}
subtrees satisfy the charge rule corresponding to their colour.  An
example is displayed in Figure~\ref{fig:example}, right. Note that $T$ itself is \emm not, one of its lower subtrees, so that there is no condition on the total charge. This is in agreement with the original definition of blossom trees in Section~\ref{sec:bip-init}.   Recall that the series $W_i$
and $B_i$ defined by~\eqref{WBdef}--\eqref{equ:blackgf}
% respectively
count blossom trees of charge $i$ rooted at a leaf and a bud, respectively.

The following lemma, dealing with blossom trees of total charge $1$,
% {immediate in the case $m=1$}
is straightforward.
% and will be% proveuseful later on:
As we shall see, these trees  are usually much simpler to handle than trees of larger charge.

\begin{Lemma}\label{lem:re-root}
  Let $T$ be a tree of total charge one, and let $e$ be an edge of
  $T$.  Then the black charge rule is satisfied at $e$ if and only if
  the white charge rule is satisfied at $e$. 

  If, in addition, $T$ is a
  blossom tree, then both charge rules are satisfied at every edge. Hence
  re-rooting~$T$ at another half-edge yields again a blossom tree.
\end{Lemma}
%
% \begin{Lemma}\label{lem:re-root}
%   Let $T$ be a tree of total charge $m\geq1$.  Assume moreover that
%   the degrees of all vertices are multiples of $m$. (Observe that this
%   is no condition if $m=1$.)
  
%   Then the black charge rule is satisfied at the edge $e$ of $T$ if
%   and only if the white charge rule is satisfied at $e$.  In
%   particular, if $T$ is a blossom tree, then both charge rules are
%   satisfied at every edge; re-rooting $T$ at another half-edge yields
%   again a blossom tree.
% \end{Lemma}

% ============================================================
\subsection{Balanced trees and their closure}
\label{subsec:closure}
%============================================================
In this section we define the \emph{closure} as a mapping $\phi$ from
certain  trees  with total charge $k\geq 1$ to $k$-leg maps.
This closure is the same as in earlier texts using the idea of
conjugacy classes of trees \cite{BMS00, BDFG02c, domi, PS02, Sch97,
gillou}.

Let $T$ be a tree with total charge $k\geq1$. Its half-edges form a
cyclic sequence around the tree in counterclockwise direction. Buds
and leaves can be matched in this cyclic sequence as if they were
respectively opening and closing brackets. More precisely, first match
the pairs made of a bud and a leaf that are immediately consecutive in
the cyclic sequence. Then, forget matched buds and leaves, and repeat
until there is no more possible matching. In view of the number of
buds and leaves in the original cyclic sequence, exactly $k$ leaves remain
unmatched. These leaves are called the \emph{single leaves} of $T$.
A tree %rooted at a leaf
is  \emph{balanced} if its root is
single; this forces the root to be a leaf.

\begin{figure}[htb]
\begin{center}
\epsfig{file=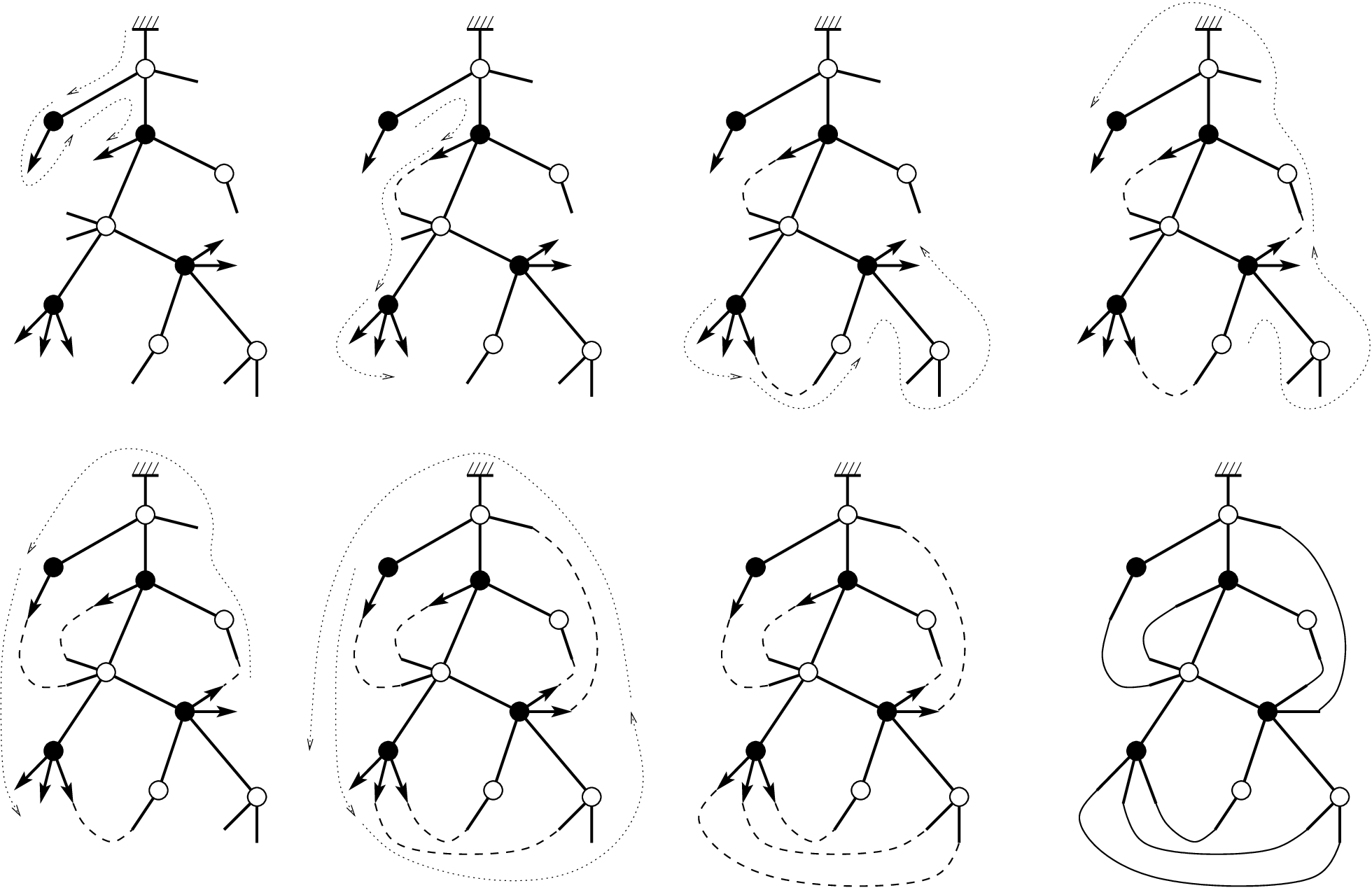, width=14cm}
\end{center}
\caption{The closure of a balanced tree having total charge $1$.}
\label{fig:closure}
\end{figure}

The \emph{closure} $\phi(T)$ of a
balanced tree~$T$ with charge $k\ge 1$ is obtained
by merging  each matched pair bud/leaf into an single edge (going from the bud to the leaf in counterclockwise direction around the
tree). The root remains unchanged, and is thus a leg of the map. An example is shown in Figure~\ref{fig:closure}. More generally, we will sometimes apply the same closure operation to a tree $T$ with total charge $k\ge 1$ in which the root half-edge is not a single leaf (and is possibly a bud): then  the resulting map  is (exceptionally) \emm not, rooted at a leg, but at the edge that was created when matching the root.

   We  have the following property.

\begin{propo}\label{pro:degrees}
Let $T$ be a balanced tree having total charge $k\ge 1$ and degree
distribution $(\lambda,\mu)$. Then $\phi(T)$ is a $k$-leg map with
degree distribution $(\lambda,\mu)$.
\end{propo}

There is an alternative description of the closure as an iterative
process, illustrated by Figure~\ref{fig:closure}: start from a
balanced tree $T$ and walk around  the infinite face in
counterclockwise order; each time a bud is immediately followed by a
leaf in the cyclic sequence of half-edges, merge them into an edge in
counterclockwise direction (this creates a new finite face that
encloses no unmatched half-edge); stop
% the course
when all buds have
been matched.  Observe that this process may in general require to
turn several times around the tree. A more efficient method is to use
a stack (that is, a first-in-last-out waiting line): store buds in the
stack as they are met and match leaves as they are met with the first
bud available in the stack. While this is algorithmically more
effective, we shall content ourselves with the previous simpler
description.

\begin{figure}[htb]
\begin{center}
\epsfig{file=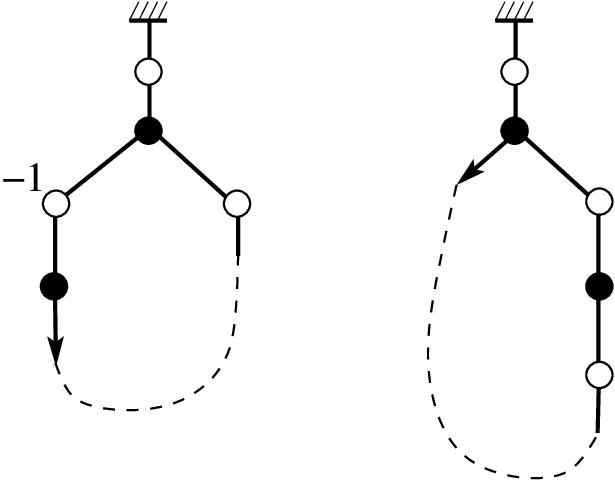, width=4.5cm}
\end{center}
\caption{Two trees that give the same $1$-leg map after closure. The left one is not a {blossom} tree, having charge $-1$ at a white-rooted subtree.}
\label{fig:non-injective}
\end{figure}

Observe that two distinct trees may yield the same map, as illustrated in Figure~\ref{fig:non-injective}. Note however that the first tree is \emm not, a blossom tree. The blossoming condition will play a key role in our construction of bijections between maps and trees.

%============================================================
\subsection{The opening of a $\boldsymbol{k}$-leg map}
\label{sec:open}
%============================================================
In this section we define the opening $\psi$ as a mapping from $k$-leg
maps to trees with total charge~$k$.
 
\begin{figure}[hbt]
\begin{center}
\epsfig{file=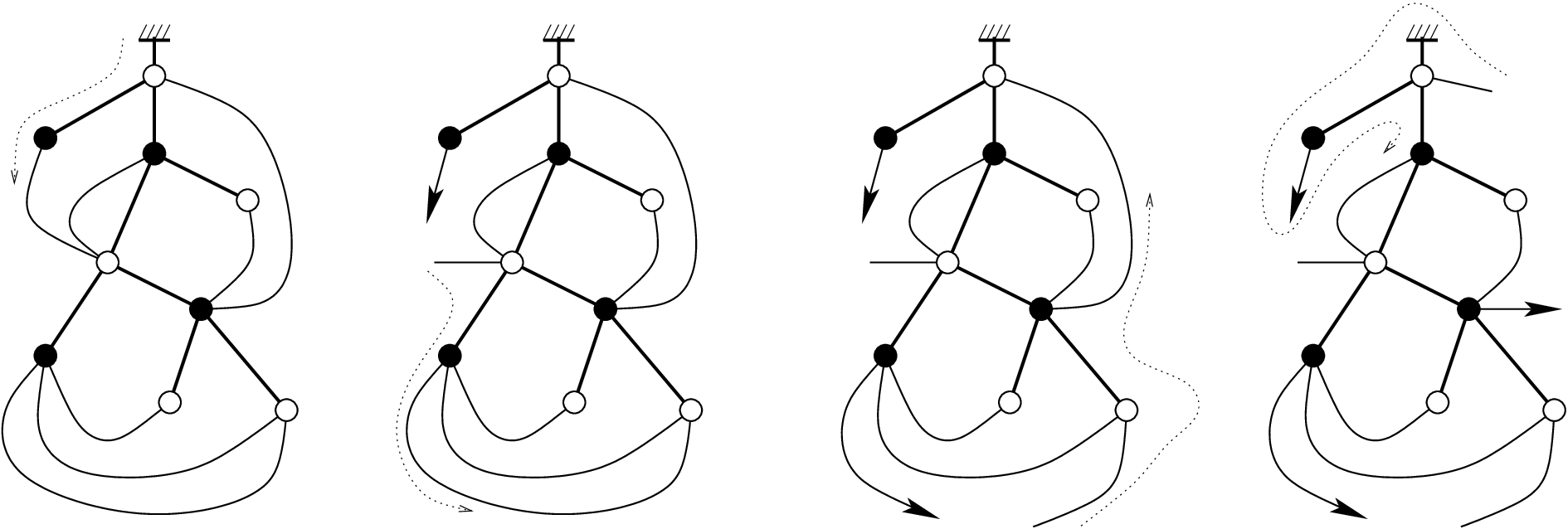, width=14cm}
\end{center}
\caption{The first few steps of the opening of a 1-leg map.
The final tree will be the one of Figure~\ref{fig:closure}.}
\label{fig:opening}
\end{figure}

The opening of a map is the result of an iterative procedure: start
from a $k$-leg map $M$ and walk around  the infinite face
in counterclockwise order, starting from the root; each time a
non-separating\footnote{An edge is \emph{separating} if its deletion
disconnects the map, \emph{non-separating} otherwise; a map is a tree
if and only if all its edges are separating.}  edge has just been
visited from its black endpoint to its white endpoint, cut this edge
into two half-edges: a bud at the black endpoint and a leaf at the
white endpoint; proceed until all edges are separating edges.  An
example is shown in Figure~\ref{fig:opening}. The opening process
cannot get stuck: as long as there are non-separating edges, a
positive even number of them are incident to the infinite face; half
of them are oriented from black to white in %the
counterclockwise
direction.  The final result is then a tree rooted at a leaf.

In view of the number of buds and leaves created, the image
$T=\psi(M)$ of a $k$-leg map is a tree with total charge $k$.
Moreover, the pairs of buds and leaves created by the
opening are in correspondence in the matching procedure of the tree described in Section~\ref{subsec:closure}.
  Hence the tree is balanced, and we have the following proposition.

\begin{propo}\label{pro:closeopen}
  The closure is inverse to the opening: for any $k$-leg map
  $M$, $\phi(\psi(M))=M$.
\end{propo}

The above proposition, combined with Proposition~\ref{pro:degrees}, guarantees that the opening~$\psi$ is a bijection between maps with $k$ legs and their image by $\psi$. The problem is now to characterize this image, which is a set of balanced trees with total charge $k$. In the next section, we characterize the image by $\psi$ of several families of maps, including those having $k=1$ or $k=2$ legs.

%%%%%%%%%%%%%%%%%%%%%%%%%%%%%%%%%%%%%%%%%%%%%%%%%%%%%%%%
\section{Bijections between maps and balanced blossom trees} 
\label{sec:bijection}
%%%%%%%%%%%%%%%%%%%%%%%%%%%%%%%%%%%%%%%%%%%%%%%%%%%%%%%%

%============================================================
\subsection{The fundamental case of one-leg maps}
%============================================================
\begin{Theorem}\label{thm:closure}
  Closure and opening are inverse bijections between balanced blossom
  trees of total charge $1$ and $1$-leg maps. Moreover, they preserve
  the degree distribution $(\lambda,\mu)$.
\end{Theorem}
In order to prove Theorem~\ref{thm:closure}, we first exhibit a
bijective decomposition of $1$-leg maps into one or two smaller
$1$-leg maps. Then, we present a related decomposition of balanced
blossom trees of total charge~$1$. We observe that the two
decompositions are \emph{isomorphic}, so that they induce a recursive
bijection between $1$-leg maps and balanced blossom trees of total
charge~$1$.
   Once the existence of a bijection has thus been established, we want to
identify it as the opening of maps. We observe that the opening
transforms the decomposition rules of maps into the decomposition
rules of trees. In view of Propositions~\ref{pro:degrees}
and~\ref{pro:closeopen}, this proves that closure and opening realize
this recursive bijection and preserve the degree distribution.

% ============================================================
\subsubsection{Decomposition of one-leg maps}
\begin{figure}[htb]
\begin{center}
\epsfig{file=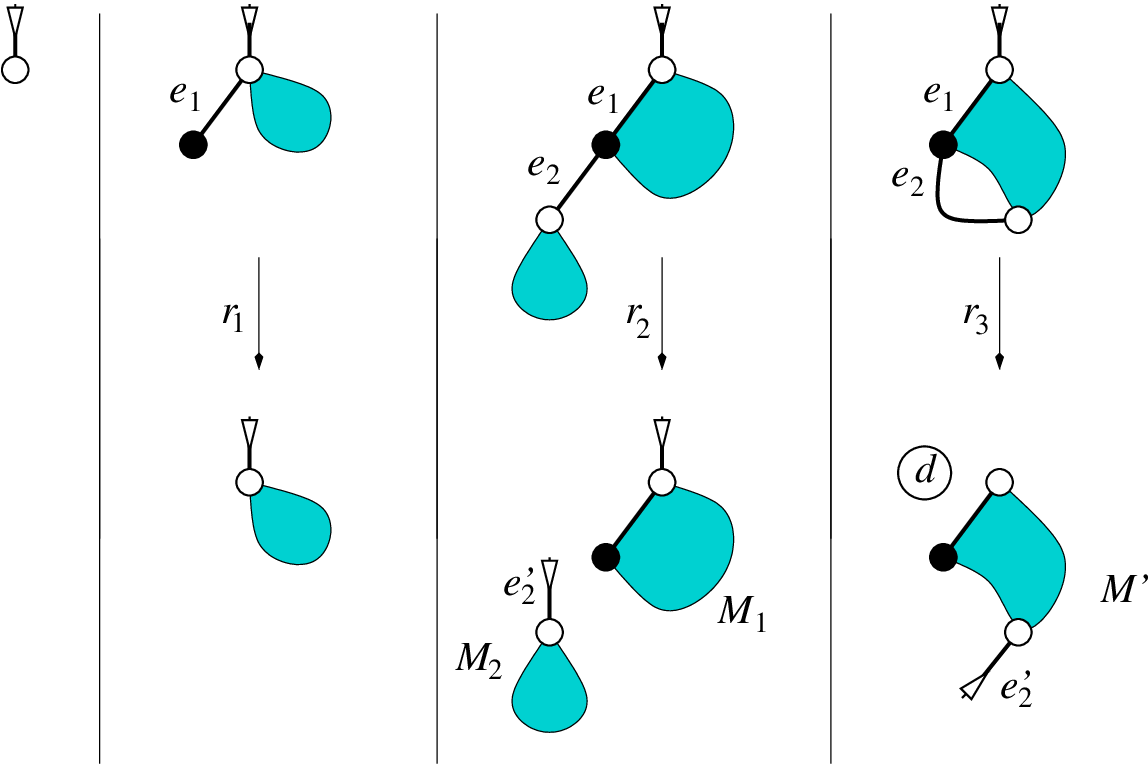, width=11.5cm}
%\rule{25mm}{0mm}
\end{center}
\caption{The decomposition of a one-leg map $M$ (shown at the top).}
\label{fig:deco-map}
\end{figure}

%============================================================
We partition the set $\mathcal M$ of one-leg maps into four disjoint
subsets and define a bijective decomposition for each of these
subsets, as illustrated by Figure~\ref{fig:deco-map}.
This is close  to the classical recursive description of planar maps by deletion of the root edge, à la Tutte. We associate two
parameters with each map $M$:
% its degree distribution $(\lambda,\mu)$,
 {the number $h(M)$ of half-edges (where an edge counts for two half-edges)}
and the \emph{reduced degree} of its infinite face, denoted by $d(M)$
(\emph{reduced} means that the root leaf is not counted, so that
$d(M)$ is even).  These parameters will be used to check that the
recursive decomposition of one-leg maps is isomorphic to that of
balanced blossom trees given in Section~\ref{sec:dec-trees}.
% further

The four subsets are defined by considering the first and second edge
$e_1$ and $e_2$ that are met when walking around the 
infinite face in counterclockwise direction, starting from the root.
We denote by $M_0$ the map reduced to a single
vertex with one leg. We denote by
${\mathcal M}^+$ the set ${\mathcal M}\setminus\{M_0\}$.

\begin{itemize}
\item \emph{$M=M_0$.}  This is the base case, where the recursive
decomposition stops. {The map~$M_0$ has one half-edge, and its infinite face has reduced degree  zero.}
  
\item \emph{$M\in\mathcal{M}_1$: $e_1$ is a dangling edge}, as defined in Section~\ref{sec:legs}.  Let
  $r_1(M)$ be the map obtained from $M$ by removing $e_1$ and its
  black endpoint.
  
  The mapping $r_1$ is a bijection from $\mathcal{M}_1$ to
$\mathcal{M}$.  Moreover {$h(M)= h(r_1(M))+2$ and} $d(M)=d(r_1(M))+2$.
  
\item \emph{$M\in\mathcal{M}_2$: $e_1$ is not dangling and $e_2$ is a
    separating edge.}  Detach $e_2$ from its black endpoint and turn
    it into a leaf $e'_2$. Two components remain: a one-leg map $M_1$
    rooted at the root leaf of $M$, and a one-leg map $M_2$ rooted at
    the new leaf $e'_2$.  Let $r_2(M)=(M_1,M_2)$.
  
  The mapping $r_2$ is a bijection from $\mathcal{M}_2$ to
  $\mathcal{M}^+ \times \mathcal{M}$.
  Moreover,  {$h(M)= h(M_1)+h(M_2)+1$ and} $ d(M)=d(M_1)+d(M_2)+2$.
  
\item \emph{$M\in\mathcal{M}_3$: $e_1$ is not dangling and $e_2$ is
    not a separating edge.}  Let $M'$ be obtained as follows: delete
    the root leaf of $M$, detach $e_2$ from its black endpoint,
    turn it into a leaf~$e'_2$, and re-root the map on this new leaf
    $e'_2$. Let $r_3(M)=(d(M),M')$.  Observe that  {$h(M)= h(M')+2$ and} $d(M')\geq d(M)$.

  The mapping $r_3$ is a bijection from $\mathcal{M}_3$ onto the set
  $\{(d,M')\mid M' \in \mathcal{M}^+, \ d \ \hbox{even}, \ 2\leq d\leq
  d(M')\}$.
 This mapping is reversible, since the value of $d$ tells
 us where to insert the new root:
 exactly $d-2$ edges after the root leaf of $M'$, in counterclockwise order around the infinite face.
\end{itemize}

\begin{figure}[h]
\begin{center}
\epsfig{file=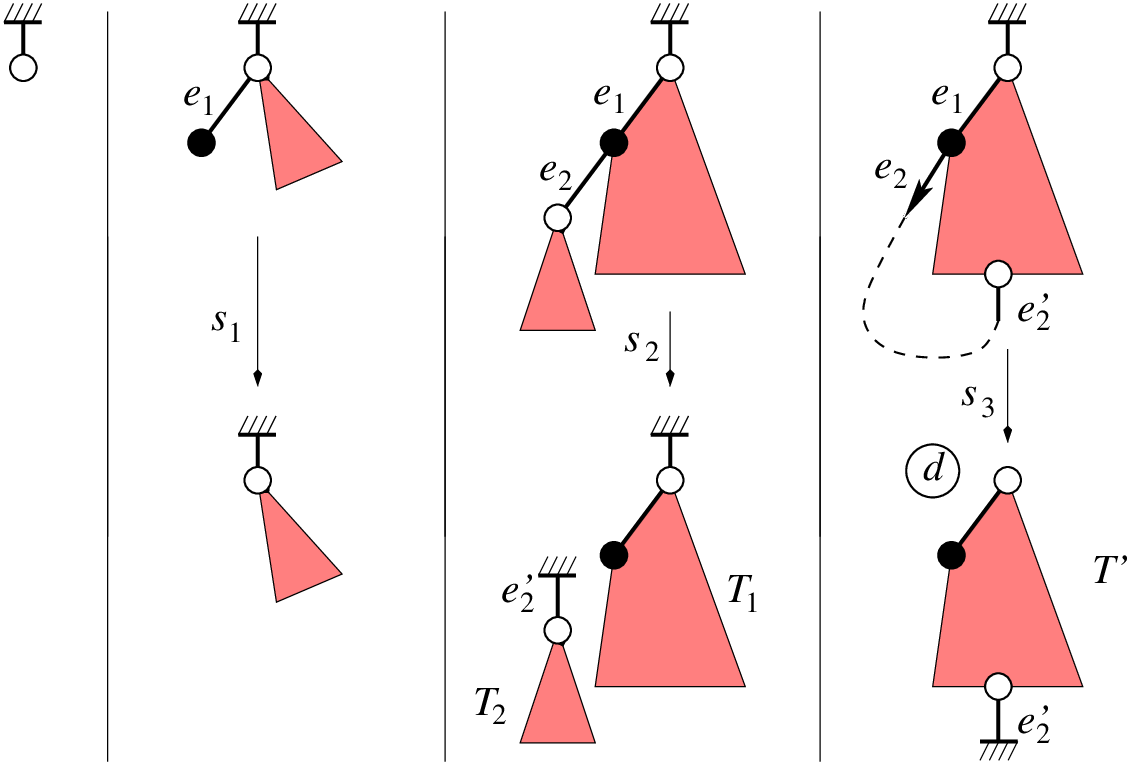, width=11.5cm}
\end{center}
\caption{The decomposition of a balanced blossom tree $T$
  of total charge~1  (shown at the top). }
\label{fig:deco-trees}
\end{figure}

% ============================================================
\subsubsection{Decomposition of balanced blossom trees with total charge one}
\label{sec:dec-trees}
%============================================================
We partition the set $\mathcal T$ of balanced blossom trees of total charge one
into four disjoint subsets and define a bijective decomposition for
each of these subsets, as illustrated by
Figure~\ref{fig:deco-trees}. The two parameters we retain for a tree
$T$ are again
% its degree distribution $(\lambda,\mu)$,
{the number of half-edges (where an edge counts for two half-edges)} and the reduced
degree of the closure of $T$, that is, $\delta(T):=d(\phi(T))$.

The four subsets of trees are defined by examining the first and
second edges or half-edges $e_1$ and $e_2$ that are met after the root leaf when walking
around the infinite face in counterclockwise direction.
% , starting from the root
Observe that $e_1$ cannot be a leaf, since $T$ is balanced and   has only one single leaf.

We denote by $T_0$ the tree reduced to a vertex with a root leaf, and
by $\mathcal{T}^+$ the set $\mathcal{T}\setminus\{T_0\}$.

\begin{itemize}
  
\item \emph{$T=T_0$.}  This is the base case, where the recursive
decomposition stops. The {tree $T_0$ has one half-edge and} the infinite face of $\phi(T_0)$ has reduced degree zero.
  
\item \emph{$T\in\mathcal{T}_1$: $e_1$ is a dangling edge.}  Let
  $s_1(T)$ be the tree obtained from $T$ by removing $e_1$ and its
  black endpoint.
  
  The mapping $s_1$ is a bijection from $\mathcal{T}_1$ to
  $\mathcal{T}$.
Moreover, {$h(T)=h(s_1(T))+2$ and}   $\delta(T)=\delta(s_1(T))+2$.

\item \emph{$T\in\mathcal{T}_2$: $e_1$ is not dangling and $e_2$ is
    not a bud.} 
Detach $e_2$ from its black endpoint
  and turn it into a leaf $e'_2$; two components remain, a tree $T_1$
  rooted at the root leaf of $T$, and a second tree $T_2$ rooted at
  the new leaf $e'_2$ (with the notation introduced
  % at the beginning of   this section,
  in Section~\ref{sec:legs}, $T_2=T^\circ_{e_2}$).
Let $s_2(T)=(T_1,T_2)$.

  We wish to prove that $T_1$ and $T_2$ are balanced blossom trees of
total charge one. Recall that $T$ is balanced, and examine the closure
of $T$.
% Since the subtree $T_2$ immediately follows the root in the infinite face of $T$,
By definition of $e_1$ and $e_2$, no bud of~$T_1$ can be matched to a leaf of
$T_2$.  Moreover, none of the (non-root) leaves of $T_2$ is single in~$T$ (because
$T$ has a unique single leaf, which is its root).  Therefore, all
leaves of~$T_2$ are matched to buds of $T_2$ in the closure of
$T$. Since $T$ is a blossom tree, the charge of $T_2 $ is nonnegative,
so that no bud of $T_2$ can match a leaf of $T_1$. In other words, the
closure takes place independently in $T_1$ and $T_2$. This proves that
$T_1$ and $T_2$ are balanced, and have total charge $1$ (and hence
charge $0$). In particular, the lower charge at $e_1$ is not modified
in the transformation, and $T_1$ and $T_2$ are blossom trees.

  Hence, the mapping $s_2$ is a bijection from $\mathcal{T}_2$ to
$\mathcal{T}^+ \times \mathcal{T}$.
Moreover, {$h(T)=h(T_1)+h(T_2)+1$ and} $\delta(T)=\delta(T_1)+\delta(T_2)+2$.
  
\item \emph{$T\in\mathcal{T}_3$: $e_1$ is not dangling and $e_2$ is a
  bud.}  Let $e'_2$ be the leaf to which $e_2$ is matched in the
  closure $\phi(T)$.  Let $s_3(T)=(\delta(T),T')$ where $T'$ is
  obtained from $T$ by deleting the root leaf and the bud $e_2$, and
  re-rooting the tree on the leaf $e'_2$.
  
  The tree $T'$ is clearly balanced, and has total charge $1$. We wish
  to prove that it satisfies the charge rules of blossom trees.
  According to Lemma~\ref{lem:re-root}, the tree obtained by
  re-rooting $T$ at $e'_2$ satisfies both charge rules at every edge,
  so that it is sufficient to consider the effect of deleting from this tree the bud
  $e_2$ and the root leaf of $T$.  Since these two deleted half-edges
  are incident to $e_1$, this edge is the only one where the charges
are modified. Let $i\geq0$ be the charge of $T^\circ_{e_1}$. Then
  $T^\bullet_{e_1}$ has charge $1-i$, and the charges of
  $T'^\circ_{e_1}$ and $T'^\bullet_{e_1}$ are respectively $i-1$ and
  $2-i$. It is thus sufficient to prove that $i\geq1$. But if $i$ was
  equal to $0$, a bud of $T^\circ_{e_1}$
  % $T_1$
  would match the root leaf of $T$, which
  would not be balanced.

  Observe that   {$h(T)= h(T')+2$ and} $\delta(T')\geq \delta(T)$. 
  The mapping $s_3$ is a bijection from $\mathcal{T}_3$ onto the set
  $\{(d,T')\mid T' \in {\mathcal  T}^+, \ d \ \hbox{even} , 2\leq d\leq
  \delta(T')\}$.    The value of $d$ tells
 us where to insert the new root:
 exactly $d-2$ edges after the root leaf of $\phi(T')$, in counterclockwise order around the infinite face.
%The value of $d$ indicates the position of the edge   $e_1$ in the infinite face of the closure $\phi(T')$.

\end{itemize}

% ============================================================
\subsubsection{Conclusion}
% ============================================================
 The decompositions of Figures~\ref{fig:closure} and~\ref{fig:opening} are ``isomorphic''. {They transform in the same way the two key parameters: the number of half-edges and the parameter $d(\cdot)$ (resp. $\delta(\cdot)$)  for maps (resp. for trees).}
  % Moreover, these figures show %that the degree distribution $(\lambda, \mu)$is altered in the same way by both decompositions.
  This implies  that there 
exists a bijection $\bar \psi$ between $1$-leg maps and balanced
blossom trees of total charge~$1$, {preserving these parameters}. %, with reverse bijection $\bar \phi$.
More precisely,  this bijection is defined by induction of the number of edges as follows:
  \begin{itemize}
  \item  $\bar \psi(M_0)=T_0$,
  \item  if $M$ is of type $1$ then    $\bar \psi(M)= s_1^{-1}\circ \bar \psi \circ r_1(M)$,
  \item  if $M$ is of type $2$ and $r_2(M)=(M_1, M_2)$ then    $\bar \psi(M)= s_2^{-1}\left(\bar \psi( M_1), \bar \psi (M_2)\right) $,
        \item  if $M$ is of type $3$ and $r_3(M)=(d,M')$ then    $\bar \psi(M)= s_3^{-1}\left(d, \bar \psi (M')\right) $.
  \end{itemize}

%, preserving the degree distribution.
Now to prove that $\psi$ and $\bar \psi$ coincide, we proceed inductively again, by checking that $\psi$ satisfies the above recursive description as well.
This is  an easy task, based on Figures~\ref{fig:deco-map} and~\ref{fig:deco-trees}. In fact, our decomposition of trees has been precisely designed to ``fit'' with the decomposition of maps via the opening procedure.
%{for instance, if a map $M$ is of the third type, the associated tree $T=\psi(M)$ is of the third type as well, and, denoting $r_2(M)=(M_1, M_2)$,  the openings of the maps $M_1$ and $M_2$ are $T_1$ and $T_2$, where $(T_1,T_2)=s_2(T)$.}
Since $\psi$ and $\bar \psi$ coincide, $\psi$ is a bijection between 1-leg maps and balanced blossom trees of total charge $1$. By Proposition~\ref{pro:degrees}, this bijection preserves the degree distribution. By Proposition~\ref{pro:closeopen}, its inverse is the closure $\phi$ of trees. This completes the proof of Theorem~\ref{thm:closure}.

%==========================================================
\subsection{Two-leg maps and other families}
%extensions
%=================================================================
  
%
%Theorem~\ref{thm:closure} extends to balanced blossom trees of total charge $2$ and $2$-leg maps.

We now present a version of Theorem~\ref{thm:closure} for   $2$-leg maps.

\begin{Theorem}\label{thm:closure-2}
  Closure and opening are inverse bijections between balanced blossom
  trees of total charge $2$ and $2$-leg maps. Moreover, they preserve
  the degree distribution $(\lambda,\mu)$.
  %There exists a bijection preserving the degree  distribution between balanced blossom trees of total charge  $2$ and $2$-leg maps. 
\end{Theorem}

\begin{proof}
  Take a balanced blossom tree $T$ of total charge $2$.
  % Observe that,
  Beyond its root leaf $r$, another leaf, say~$e_1$, is single. It ends up in the infinite face of the map $\phi(T)$ when closing the edges of~$T$.
  The idea is to replace one of the single leaves of $T$ by a black vertex of degree $1$, so as to obtain a tree with total charge $1$, and then apply Theorem~\ref{thm:closure}. However, transforming the leaf $e_1$ in this way may yield a tree that is not a blossom tree, as can be easily seen. What we do instead is to
  replace \emm the root leaf, $r$ by
  % an edge $\bar r$ ending with a
  a black vertex of degree $1$ -- more precisely, by a dangling edge $e_0$ ending at
  marked black vertex -- and
  choose $e_1$ as the new root. We obtain a tree $T'$ of total charge $1$, rooted at a leaf and carrying a marked black vertex of degree $1$.

  We claim that $T'$ is a  blossom tree. First, the charge condition is clearly satisfied at the new edge $e_0$. Then, for any edge $e$ of $T'$ that was already in $T$, one of the subtrees $T'^\circ_e$ and~$T'^\bullet_e$ is a lower subtree in $T$, and thus satisfies the charge condition. By the first part of Lemma~\ref{lem:re-root}, this implies that $T'$ is a blossom tree.  Moreover, the pairs of buds and leaves that
are matched by closure are the same in $T'$ as in $T$ (the closure
does not depend on the root). In particular~$e_0$ and $e_1$
{remain} in the infinite face in $\phi(T')$, and $T'$ is balanced.

We claim that the transformation $f$ that transforms~$T$ into $T'$ is a bijection between balanced blossom trees of total charge $2$ and
balanced blossom trees $T'$ of total charge $1$ with a marked black vertex
of degree~$1$ \emph{that ends up in the infinite face of $\phi(T')$
  when the edges of $T'$ are closed}: indeed, for any tree $T'$
satisfying the latter properties, one obtains a balanced blossom tree
$T$ of total charge $2$ by replacing the marked black vertex by a leaf
which becomes the new root; the proof  that $T$ is a blossom
tree is similar to the proof that $T'$ was a blossom tree above.

 Now let us apply Theorem~\ref{thm:closure}: the closing procedure
$\phi$ induces a bijection between the trees~$T'$ satisfying
the above properties and 1-leg maps $M'$ having a marked
black vertex of degree $1$ in the infinite face. In $M'$, replace the
marked vertex by a leaf that becomes the new root: the resulting map
$M$ is precisely $\phi(T)$, where $T=f^{-1}(T')$, and $M$
has  the same
degree distribution as~$T$. This concludes the proof.
\end{proof}

The case of $k$-leg maps with $k\ge 3$ is more complex, but the
following result will be sufficient for our purpose.

\begin{Theorem}\label{thm:closure-m}
  Let $m\geq 3$ be an integer. Closure and opening are inverse
  bijections between balanced blossom trees of total charge $m$ in
  which the degrees of vertices are %all
  multiples of $m$, and $m$-leg
  maps satisfying the same condition.  %Moreover
  The degree   distribution $(\lambda,\mu)$ is preserved.
\end{Theorem}

The proof of this theorem is based on the following  simple lemma. 

\begin{Lemma}\label{lem:multipleofm}
Let $T$ be a  tree in which the degrees of all vertices
  are multiples of $m$. Then $m$  divides the total charge of $T$.

In particular, a lower subtree with a white root vertex satisfies the
charge rule if and only if its  charge is at least $m-1$ (a stronger
condition).  A lower subtree with a black root vertex satisfies the
charge rule if and only if its total charge is at most $m$ (a weaker condition).

%Let $T$ be a blossom tree in which the degrees of all vertices are  multiples of $m$. Let~$T'$ be obtained from $T$ by replacing up to  $m-1$ non-root leaves of $T$ by edges ending with a black vertex.  Then $T'$ is still a blossom tree.
\end{Lemma}

\begin{proof} Let $T$ be a  tree in which the degrees of all vertices
are multiples of $m$. Let $(\lambda,\mu)$ denote the vertex
distribution of $T$. Observe that $m$ divides $|\lambda|$ and
$|\mu|$. Then $m$ also divides the total charge of $T$, because the
latter is the difference between the numbers $n_\ell$ of leaves and~$n_b$ of
buds, which satisfy $n_b+|\lambda|=n_\ell+|\mu|$ (as observed in Section~\ref{sec:legs}).

Now the vertices of any lower subtree of $T$ have a degree that
  is a multiple of $m$. Thus if $T^\circ$ is  a lower subtree of $T$ with a
  white root vertex that satisfies the charge rule, the total charge
  of $T^\circ$ is at least $1$, hence at least $m$, hence the charge  of
  $T^\circ$ is   in fact 
  at least $m-1$.  If $T^\bullet$ is a lower subtree of $T$ with a
  black root vertex and a charge at most $m$, its total charge is at
  most $m-1$, hence at most $0$ by the first part of the lemma; this
  means that the charge of $T^\bullet$ is at most $1$, so that the charge rule is satisfied.
%Assume moreover that $T$ is a blossom tree. Then any lower subtree$T^\circ_e$ of $T$ has by definition a charge $c_\circ\geq0$.According to the previous discussion the total charge $c_\circ+1$ of$T^\circ_e$ is divisible by $m$ so that in fact $c_\circ\geq m-1$.Since replacing a leaf by an edge ending with a black vertex decreasesthe charge by one, up to $m-1$ leaves can be replaced without any riskto violate a white charge rule.
\end{proof}

\begin{proof}[Proof of Theorem~\ref{thm:closure-m}.] Let  $T$ be 
a balanced blossom tree of total charge $m$ in which the degrees
of vertices are multiples of $m$.
% Observe that
Its $m$ single
leaves (among which the root) end up in the infinite face of the map
$\phi(T)$ when closing the edges of $T$. Replace each of the single
leaves of $T$, except the root, by a {dangling} edge ending with a black vertex
of degree~$1$. This gives a balanced
% blossom
tree~$T'$ of total charge~$1$, in which all vertex degrees are
multiples of $m$ apart from the $m-1$ %marked
black
vertices of degree~$1$; moreover, these vertices end up in the infinite face of
$\phi(T')$ when we close the edges of $T'$.

We claim that  $T'$ is a blossom tree: given that the charge of any lower
  subtree has decreased (by at most $m-1$) during the transformation $T\mapsto T'$, there is obviously no problem for
  black rooted subtrees. But there is no problem either for  white rooted subtrees, whose
  charge in~$T$ was at least $m-1$ by  Lemma~\ref{lem:multipleofm}, and thus remains non-negative.
Conversely, starting from a balanced blossom tree $T'$ satisfying all the
above conditions, replacing each marked black vertex of degree $1$ by
a leaf yields a balanced tree $T$ of total charge $m$ in which all
degrees are multiples of $m$, and this tree is a blossom tree (the
argument is similar to the above one, but we now pay attention to
black rooted subtrees). Hence the transformation
$f:T \mapsto T'$ is bijective between these two classes of
trees. By Theorem~\ref{thm:closure}, the closure $\phi$ sends
bijectively the trees $T'$ to 1-leg maps $M'$ having all vertices of
degree a multiple of $m$, apart from $m-1$ black vertices of degree
$1$ in the infinite face.  Converting
each of these into a leg gives the map $\phi(T)$, where
$T=f^{-1}(T')$, and this concludes the proof of  Theorem~\ref{thm:closure-m}.
\end{proof}

The following interesting variation was observed, in the non-bipartite
case, by by Bouttier \emph{et
al.}~\cite{BDFG02c}.   It describes a class of maps that are in
bijection with trees \emph{not subject to any balance condition},
which are much easier to count.

\begin{Theorem}
  There is a one-to-one correspondence between (not necessarily
  balanced) blossom trees of charge $1$ rooted at a leaf and $1$-leg
  maps having a marked black vertex of degree~$1$, that sends a tree with degree distribution $(\lambda, \mu)$ to a map with degree distribution $(\lambda, \mu')$, with $\mu'= 1\cdot \mu$.
  
  Similarly, there is a one-to-one correspondence between blossom
  trees of charge $1$ rooted at a bud and $1$-leg maps having a marked
  white vertex of degree $1$, that sends a tree with degree distribution $(\lambda, \mu)$ to a map with degree distribution $(\lambda', \mu)$, with $\lambda'= 1\cdot \lambda$.
\end{Theorem}

\begin{proof}
Consider a blossom tree $T$ of charge $1$, rooted at a leaf, not necessarily balanced. As in the proof of Theorem~\ref{thm:closure-2}, replace the
root leaf by a marked black vertex of degree $1$ and re-root the
resulting tree $T'$ of charge one on its unique single leaf. The argument
used for Theorem~\ref{thm:closure-2} shows that this transformation is a bijection between blossom trees of charge $1$ rooted at a leaf and balanced blossom trees with a marked black vertex of degree $1$. One then applies the bijection of Theorem~\ref{thm:closure}.   A similar argument proves the second statement, upon creating a  white vertex of degree $1$ at the root bud. 
\end{proof}

%%%%%%%%%%%%%%%%%%%%%%%%%%%%%%%%%%%%%%%%%%%%%%%%%%%%%%%%
\section{Counting balanced trees}
\label{sec:conjugue}
%%%%%%%%%%%%%%%%%%%%%%%%%%%%%%%%%%%%%%%%%%%%%%%%%%%%%%%%
In view of Eq.~\eqref{equ:maps-pattes} and Theorems~\ref{thm:closure}
and~\ref{thm:closure-2}, our objective is now to count balanced
blossom trees of total charge $1$ or $2$. In this section, we
establish bijections that allow us to express their \gfs\ in terms of
the \gfs\ $W_i$ and $B_i$ of (not necessarily balanced) blossom
trees, defined by~\eqref{WBdef}--\eqref{equ:blackgf}. These bijections are adaptations to the bipartite case of the
bijections presented in~\cite{BDFG02c}.

For $i\in\zs$, let $\mathcal{W}_i$ be the set of blossom trees of
charge $i$ rooted at a leaf (so that the total charge is $i+1$).
Similarly, let $\mathcal{B}_i$ be the set of blossom
trees of charge $i$ rooted at a bud (their total charge is $i-1$).
These trees are respectively counted by the series $W_i$ and $B_i$ of
Eqs.~\eqref{WBdef}--\eqref{equ:blackgf}.
Finally, for $i\geq1$, let $\mathcal{W}_i^*$ be the subset of
$\mathcal{W}_i$ formed of balanced trees.
Given  Theorems~\ref{thm:closure}
and~\ref{thm:closure-2}, we are especially interested in the enumeration of the trees of
$\mathcal{W}_0^*$ and $\mathcal{W}_1^*$. 
The following theorem implies our
Theorem~\ref{thm:main} on the enumeration of bipartite maps.
\begin{Theorem}\label{thm:bijection}
  There exists a simple degree-preserving bijection between the sets $
  \mathcal{W}_0$ and $ \mathcal{W}^*_0\cup \mathcal{B}_2.  $ There
  also exists a degree-preserving bijection between the sets $
  \mathcal{W}_1$ and $ \mathcal{W}^*_1 \cup
  \mathcal{B}_3\cup(\mathcal{B}_2)^2.  $
\end{Theorem}

Before proving this theorem, let us state a similar result which,  thanks to Theorem~\ref{thm:closure-m},  will
allow us to count certain $m$-leg maps, with $m \ge1$. It generalizes
the results obtained in~\cite{BMS00} for some bipartite maps called
\emm constellations,.

\begin{Theorem}\label{thm:constellationplusplusbifluore}
Let $m \ge1$, and let $\widehat{\mathcal{W}}_i$,
$\widehat{\mathcal{W}}^*_i$ and $\widehat{\mathcal{B}}_i$ be the
restriction of the sets $\mathcal{W}_i$, $\mathcal{W}^*_i$ and~$\mathcal{B}_i$ to trees in which the degrees of all vertices are
multiples of $m$.  There exists a simple degree-preserving bijection
between $ \widehat{\mathcal{W}}_{m-1} $ and
$\widehat{\mathcal{W}}^*_{m-1}\cup\widehat{\mathcal{B}}_{m+1}.$
\end{Theorem}
Observe that the case $m=1$ is just the first part of Theorem~\ref{thm:bijection}.

  % =====================================================
  \subsection{Trees of total charge 1: Proof of Theorem~\ref{thm:bijection} (first part)}
  % =====================================================

  The first bijection of Theorem~\ref{thm:bijection} is
extremely simple to describe: a tree $T$ belonging to $\mathcal{W}_0$
is either balanced (that is, belongs to $\mathcal{W}_0^*$) or its root
leaf is matched to a bud.
In this case, let us re-root $T$  at this bud; the resulting tree
$T'$ has still total charge~$1$ and Lemma~\ref{lem:re-root} proves
that it is a blossom tree, hence an element of $\mathcal{B}_2$.
Conversely, if we take a tree~$T'$
in $\mathcal{B}_2$ and re-root it at the leaf matched to the root bud,
we obtain a tree of~$\mathcal{W}_0$ that is not balanced.
\qed%\end{proof}

% =====================================================
\subsection{The core of a blossom tree}
% =====================================================
We would like to have a similar argument for balanced blossom trees of
total charge $2$, and, why not, of total charge $k\ge 2$. What prevents us
from doing so? Let us take a tree $T$ in~$\mathcal{W}_1$: either it is
balanced (that is, belongs to $\mathcal{W}_1^*$), or its root leaf is
matched to a bud. In this case, let us re-root $T$ at this bud to
obtain another tree $T'$: then $T'$ has still total charge $2$. Does~$T'$ satisfy the charge conditions? Well, not always... The
following lemma, whose proof is self-contained, describes in detail what may happen. It states in particular that in a blossom tree of positive total charge, the white charge condition holds at every edge.

\begin{Lemma}\label{lem:weakedges}
Let $T$ be a blossom tree with total charge $k\geq1$ and let $e$ be an edge
of $T$. Denote by $c_\bullet$ the charge of $T^\bullet_e$ and by $c_\circ$
the charge of $T^\circ_e$, so that $c_\bullet+c_\circ=k$.
\begin{itemize}
\item If the root of $T$ belongs to $T^\circ_e$ then the black charge
  rule is satisfied at $e$ (that is, $c_\bullet\leq1$) and the charge
  of $T^\circ_e$ is at least $k-1$: $c_\circ\geq k-1$.  Hence the
  white charge rule is satisfied at $e$ as well.
\item If the root of $T$ belongs to $T^\bullet_e$ then the white
  charge rule is satisfied at $e$ (that is, $c_\circ\geq0$) and the
  charge of $T^\bullet_e$ is at most $k$: $c_\bullet\leq k$.  Hence
  the black charge rule may not be satisfied at $e$ if $k\ge 2$.
\end{itemize}
If $c_\bullet\geq2$ the edge $e$ is called \emph{weak}. Weak edges may only exist if $k\ge 2$.
\end{Lemma}
The above lemma implies that if we re-root a blossom tree ``on the wrong
side'' of a weak edge, the result will not be a blossom tree.
% In other words, we have the following lemma.
This can be reformulated as follows.

\begin{Lemma}\label{lem:core}
  Let $T$ be a blossom tree with total charge $k\ge 1$, and let $T'$ be obtained by re-rooting~$T$ at a half-edge $e$. Then $T'$ is a blossom tree if and only if
there is no weak edge on the path from the root of $T$ to~$e$. In particular, if $k=1$ then $T'$ is always a blossom tree, as stated in Lemma~\ref{lem:re-root}.
\end{Lemma}

%  \noindent  {\bf C. The core of a blossom tree.} 
This leads to the notion of \emph{core}, introduced in~\cite{BDFG02c}:
let $T$ be a blossom tree with total charge $k\ge 1$, and let $S$ be the set of weak edges of $T$.
The core of~$T$ is the connected component of $T\setminus S$
containing the root.  Lemma~\ref{lem:core} can be reformulated by saying
that re-rooting a blossom tree of positive total charge gives a blossom tree if and only if the
new root is  chosen in the core. In this case, the core of the new
tree coincides with that of the original one.

Lemma~\ref{lem:weakedges} implies that  a blossom tree of  total charge $k\ge 1$ satisfies the white charge rule at every edge. Lemma~\ref{lem:re-root} implies that the converse is true when $k=1$. The following lemma relies on  the notion of core to state some kind of reciprocal when $k=2$. This will be useful in the proof of the second part of Theorem~\ref{thm:bijection}.

\begin{Lemma}\label{lem:white-charge}
Let $T$ be a tree of total charge $2$ such that the white charge condition holds at every edge. Then there exists a leaf $\ell$ of $T$ such that re-rooting    $T$ at $\ell$ gives a blossom tree.
\end{Lemma}
\begin{proof}
We first generalize the notion of weak edges and core to trees $T$ satisfying the conditions of the lemma.  We say  that an edge $e$ of $T$ is \emm weak, if it has charge $2$ at its black endpoint (and thus charge $0$ at the other end). For all non-weak edges, both charge rules hold. Let $S$ be the set of weak edges of $T$. We claim that, when deleting all edges of $S$, we obtain a collection of trees such that all trees but one have total charge $0$; the last tree has total charge $2$.  This can be proved by induction of the number of weak edges:
  let~$e$ be a weak edge of $T$ such that there is no further weak edge in~$T^\circ _e$. Then deleting $e$ leaves a component of total charge $0$ (this is~$T^\circ _e$ without its root leaf), and a tree of total charge $2$ (this is~$T^\bullet_e$ without its root bud)  whose weak edges are those of $T$, except $e$. We then repeat the procedure to conclude. We define  the \emm core, of $T$ as the unique component of $T\setminus S$ with total charge~$2$.

Observe that if $T$ is also a blossom tree, the definition of weak edges coincides with that of Lemma~\ref{lem:weakedges}. Moreover, if $e$ is a weak edge then the root of $T$ belongs to $T^\bullet_e$ by Lemma~\ref{lem:weakedges}. This shows that both notions of core coincide as well.

Let us now return to a tree $T$ that satisfies the white charge rule at every edge but is not necessarily a blossom tree. Let $\ell$ be one of the leaves of the core (there are at least $2$ of them), and let $T'$ be obtained by re-rooting $T$ at $\ell$. We claim that $T'$ is a blossom tree. Since the white charge rule is satisfied at every edge of $T$ (hence of $T'$), and the black charge rule at every non-weak edge, we just need to consider weak edges, and to prove that the lower subtree of $T'$ at a weak edge $e$ is always $T^\circ_e$. If this were not the case, deleting $e$ would leave a component of total charge $2$ not containing the core, which is impossible. This concludes the proof.
\end{proof}

%\smallskip
%The above discussion shows that it is difficult to characterize thetrees obtained by re-rooting an unbalanced blossom tree $T$ of totalcharge $k\ge 2$ at the bud that matches its root: if $T$ has anon-trivial core, we may obtain a tree that is not a blossom tree.

%==============================================================
\subsection  {When all degrees are multiples of
  % the total charge
  $\boldsymbol m$:
  proof of Theorem~\ref{thm:constellationplusplusbifluore}}
% ==============================================================
As already exemplified by the proof of
Lemma~\ref{lem:multipleofm}, the assumptions of
Theorem~\ref{thm:constellationplusplusbifluore} make it more difficult
to violate the charge rule.
Let us  prove that, under these assumptions, a blossom
tree $T$ with total charge $m$ has {$T$ as its  core}.  
Assume that $T$ contains a weak edge~$e$. This is only possible if $m\ge 2$. Then~$T^\circ_e$ has a
(nonnegative) charge $c_\circ=m-c_\bullet\leq m-2$ (see
Lemma~\ref{lem:weakedges}).  Equivalently,~$T^\circ_e$ is a blossom
tree with total charge $c=c_\circ+1$ in the interval $[1,m-1]$. This is impossible by Lemma~\ref{lem:multipleofm}, since all vertex degrees of $T^\circ_e$ are multiples of $m$.
% Let$n_b$ be the number of buds of~$T^\circ_e$, $n_\ell=n_b+c$ its numberof leaves and $(\lambda,\mu)$ its degree distribution. Then$n_b+|\lambda|=n_b+c+|\mu|$. But $c$ belongs to $[1,m-1]$, and thiscontradicts the assumption that $m$ divides both $|\lambda|$ and$|\mu|$.
Hence $T$ has no weak edge. 

We can now prove Theorem~\ref{thm:constellationplusplusbifluore} by adapting the argument  used above
  % in the case $k=1$
  to prove the first part of Theorem~\ref{thm:bijection}. A tree belonging to  $\widehat{\mathcal{W}}_{m-1} $ is either balanced (that is, belongs to $\widehat{\mathcal{W}}^*_{m-1}$) or its root leaf is matched to a bud. In this case, we re-root $T$ at this bud. Since $T$ has no weak edge, the resulting tree~$T'$ is a blossom tree, of total charge $m$, hence an element of 
$\widehat{\mathcal{B}}_{m+1}$. Conversely, if we take a tree in $\widehat{\mathcal{B}}_{m+1}$ and re-root it at the leaf matched to the root bud, we obtain an unbalanced tree of~$\widehat{\mathcal{W}}_{m-1} $.
This concludes the proof of Theorem~\ref{thm:constellationplusplusbifluore}.
\qed

% ==============================================================
\subsection{Trees of total charge 2: proof of Theorem~\ref{thm:bijection} (second part)}
\begin{figure}[htb]
\begin{center}
\epsfig{file=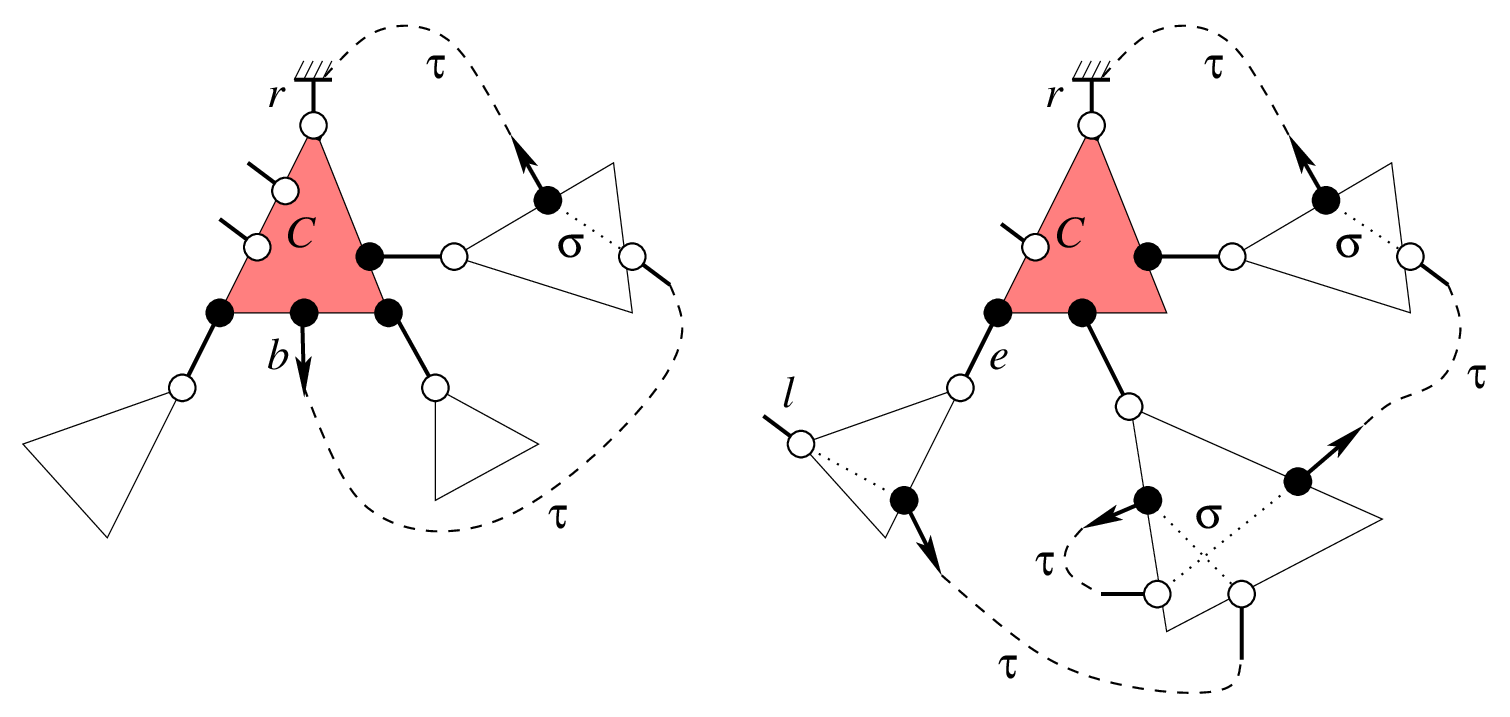, width=13cm}
\end{center}
\caption{Chains in the graphs of $\sigma$ and $\tau$, starting from the root leaf.  In the first
tree the root is sent {ultimately} to a bud $b$ of the core. In the second one the root
is sent to a single leaf $\ell$ which is not in the core.}
\label{fig:journeys}
\end{figure}

% ==============================================================
We are thus left with the  second statement of
Theorem~\ref{thm:bijection}, definitely the most difficult in this section. 
It admits (at least) two proofs, or maybe two variants of the same
proof.  First, one can adapt the not-so-easy arguments
of~\cite[Section 3.4 and Appendix B]{BDFG02c}. We choose instead to
present an explicit bijection $\rho$ from $\mathcal{W}_1\setminus \mathcal{W}^*_1$ to $\mathcal{B}_3\cup(\mathcal{B}_2)^2$, which can be roughly described
% summarized
as follows.
Let us take a blossom tree~$T$ in $\mathcal{W}_1$ having the
misfortune of being unbalanced.
%We want to associate with it anelement $\rho(T)$ belonging to $\mathcal{B}_3\cup(\mathcal{B}_2)^2$.
The tree~$T$ is rooted at a leaf and, following our earlier principle,  we want to re-root it at a bud  --- but this bud should be 
in the core, so as to obtain a blossom tree. Problems arise when the bud matched to the root
% by the closure
is not in the core.
The idea is then to go {across} the subtrees that dangle from the
 core to look further away for our dream bud.  This idea is
 schematized by Figure~\ref{fig:journeys}, and will be explained in
 greater detail below.
  Sometimes we shall reach a bud of the core, and obtain  a tree of~$\mathcal{B}_3$ by re-rooting.  Sometimes we will fail to return to the core, and then $\rho$ will return a pair of trees of
% this will yield the term
$\mathcal{B}_2$.

In order to make the construction more precise, some notations are
useful.   For every
tree $T_0$ with charge zero, let us fix an arbitrary bijection $\sigma_{T_0}$
from its buds to its leaves, ignoring the root half-edge.
The following construction is illustrated in Figure~\ref{fig:journeys} on a leaf-rooted tree.  Take a blossom tree~$T$  with total charge $2$.
%
%and assume
%{By Lemma~\ref{lem:weakedges}}, it satisfies the white charge rule at every edge.
% (this is the case for instance if $T$ is a blossom tree, .
Let $C$ be its core.
% (which can be defined, as observed above, without any reference to the root).
Then $C$ contains the root half-edge of~$T$. Let
$E$ be the set of weak edges incident to $C$. Each edge $e$ of $E$
defines a subtree $T_e^\circ$ of charge zero, rooted at a leaf. Any vertex of $T$ is
either in the core or in one of the trees $T_e^\circ$, for $e \in E$.
The bijections $\sigma_{T_e^\circ}$, for $e\in E$, induce a
bijection $\sigma$  between the buds and leaves  that are not in the core.
The {matching procedure} %closure
on~$T$ induces another bijection $\tau$ from its
non-single leaves to its buds. Observe that $\sigma$ and $\tau$ do not
depend on which half-edge (of the core) is the root of $T$.  
The graph of the bijection $\sigma$ (resp.~$\tau$) is made of oriented
edges going from buds to leaves (resp. from leaves to buds). Hence the
union of these two graphs consists of
% is made of
oriented cycles and of  chains that alternate between
buds and leaves of $T$.
Moreover, the starting points of these chains are all leaves of the core, because all other half-edges of $T$ are in the image of $\tau$ or $\sigma$. The endpoints of the chains are all buds of the core and all single leaves. The only chains that end at a (single) leaf of the core are reduced to this leaf and have length $0$.

 Now take $T$ in
$\mathcal{W}_1\setminus\mathcal{W}_1^*$ (Figure~\ref{fig:journeys}). Then $T$  is a blossom tree of total charge $2$, so that the above construction applies. The root $r$ of $T$ is a leaf of the core, hence the starting point of a chain. Since $T$ is not balanced, this chain has length $\ge 1$.
%which is not in the image of $\sigma$ ({since $r$}
% it
%is in the core) but is in the domain of $\tau$. It is thus the %origin
%starting point of a chain of length $\ge 1$.
%What is the endpoint of this chain? {As discussed above,}
% There are two cases.
 The endpoint of this chain is either 
% in the image of~$\tau$ but not in the domain of $\sigma$; that is, it is
a bud $b$ of the core (Figure~\ref{fig:journeys}, left),
or
% the endpoint is in the image of $\sigma$ but not in the domain of $\tau$; that is, it is
a single leaf $\ell$ of $T$ that is not in the core
(Figure~\ref{fig:journeys}, right).   We shall define $\rho(T)$
separately in these two cases. Let $\mathcal{W}_b$ be the subset of
$\mathcal{W}_1\setminus\mathcal{W}^*_1$ consisting of trees such that the chain
starting at 
the root ends at a bud, and let $\mathcal{W}_\ell$ be the 
complement subset.

\smallskip
\paragraph{\bf First case:}
% assume that
$T$ belongs to $\mathcal{W}_b$. 
Re-rooting the tree $T$ %on
at the endpoint $b$ yields a tree $T'$.  Since~$b$ is a bud of the core, $T'$ belongs to $\mathcal{B}_3$, and has the same core as $T$.  Set
$\rho(T)=T'$.  Conversely, take a tree $T'$ of $\mathcal{B}_3$ (hence of total charge $2$). Define
the bijections $\sigma$ and $\tau$ as above.  The union of their
graphs defines cycles and chains.  The root of $T'$ is a bud of the
core, hence the endpoint of a chain. As observed earlier, the %origin
starting point of this chain is a
leaf~$r$ of the core. Moreover, this leaf is not single, otherwise the chain would contain no bud.
% : indeed, these leaves are the only half-edgesto be in the domain of one bijection and not in the image of theother.
Re-rooting the tree $T'$ on $r$ yields a tree $T$ of
$\mathcal{W}_b$. Since the bijections
% definitions of
$\sigma$ and $\tau$ do not
depend on the root, $T$ is the only tree such that $\rho(T)=T'$.
Hence~$\rho$ is one-to-one between $\mathcal{W}_b$ and $\mathcal{B}_3$.

\begin{figure}[htb]
\begin{center}
\epsfig{file=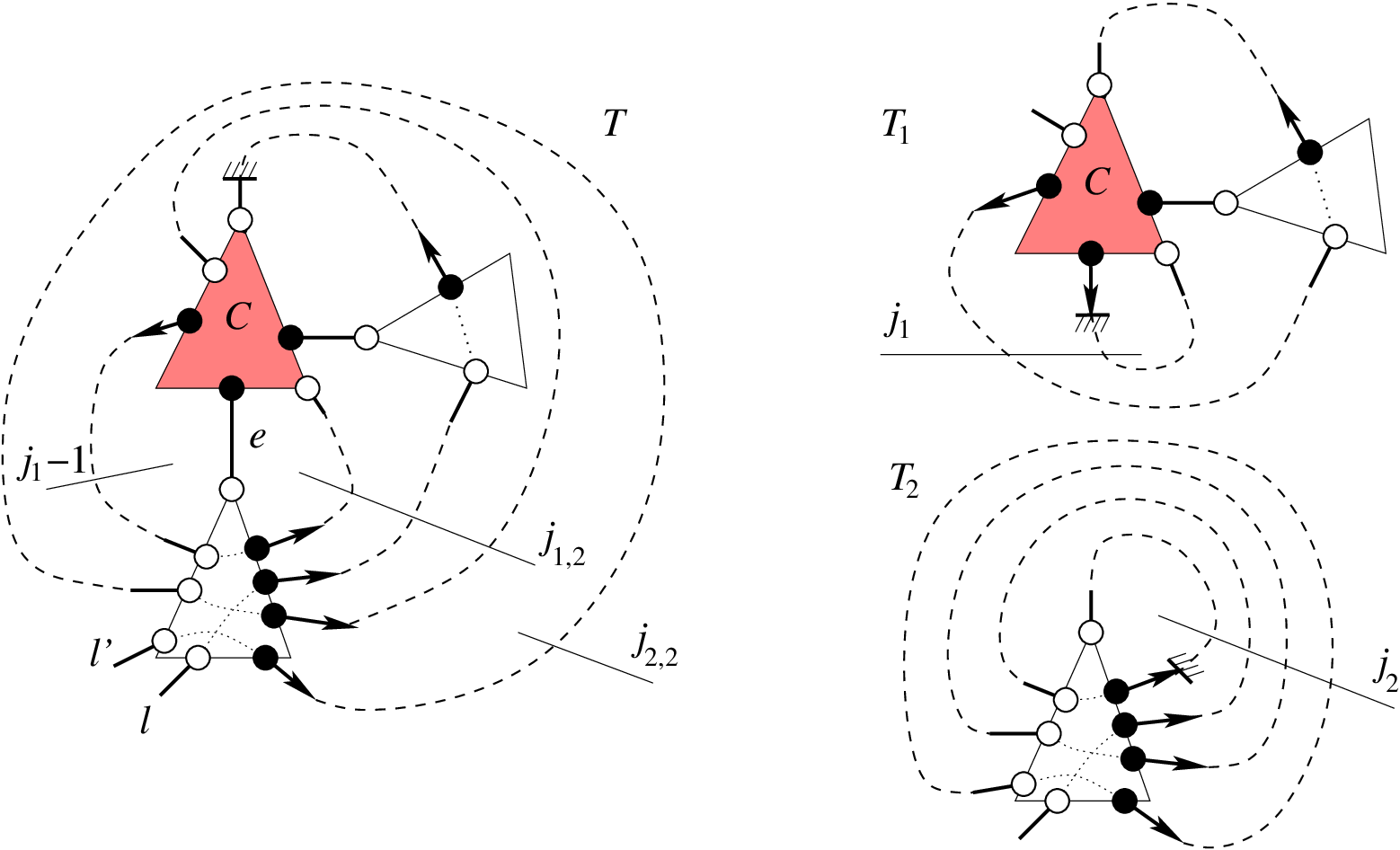, width=13cm}
\end{center}
\caption{The decomposition of a tree $T$ of $\mathcal{W}_\ell$. 
In this example, the triple $(\ell',\ell,e)$ 
 appears in counterclockwise order around $T_e^\circ$, so that
 $\rho(T)=(T_2,T_1)$.} 
\label{fig:appendix}
\end{figure}

\smallskip
 \paragraph{\bf Second case:}
% assume that
 $T$ belongs to $\mathcal{W}_\ell$. We would like to transform $T$ into a pair of trees of $\mathcal B_2$.
 %by a bijection $\bar \rho$.
 Let~$\ell$ be the single leaf ending the chain that starts at the
root $r$ of $T$.  An example is given in Figure~\ref{fig:appendix}.  The
leaf $\ell$ is not in the core.  Therefore there is a (unique)
edge~$e$ of $E$ such that~$\ell$ belongs to~$T_e^\circ$.
Let  $T_1:=T_e^\bullet$ (which is rooted at a bud) and let $T_1'$ be obtained by re-rooting $T_1$ at the leaf $r$. Then~$T_1'$ has charge $1$, and it is a blossom tree because there is no weak edge between the root and the bud that has replaced $e$. Hence by Lemma~\ref{lem:core}, the tree $T_1$ belongs to $\mathcal{B}_2$.  Consider now~$T_e^\circ$. It is a blossom tree of  charge $0$.
Let us prove that is it not balanced.
The $\sigma$-$\tau$ chain that starts from the root of $T$ and ends at the leaf $\ell$
enters~$T_e^\circ$ for the first time at a bud. This bud is matched in
$T$ to a leaf of $T_e^\bullet$, and this prevents $T_e^\circ$ from
being balanced. Hence, re-rooting~$T_e^\circ$ on the bud matched to
the root of $T_e^\circ$, as in the proof of the first part of Theorem~\ref{thm:bijection}, yields a tree of $\mathcal{B}_2$, which we denote by $T_2$.
The pair $ \rho(T)$ will be either $(T_1, T_2)$ or $(T_2, T_1)$, as we now explain.

Let us consider the closure of $T$ (Figure~\ref{fig:appendix}). 
Let
% us denote by
$j_1-1$ be the number of buds of~$T_e^\bullet$ that are
matched to a leaf of $T_e^\circ$. When walking around $T_e^\circ$ in counterclockwise order, starting from its root, one meets these leaves before reaching the single leaf $\ell$. Observe that, in the closure of~$T_e^\bullet=T_1$, the (black) root vertex then
% of $T_e^\bullet$
lies at 
 {\em depth}  $j_1$ (Figure~\ref{fig:appendix}, top right).
Let  $j_{1,2}$ be 
 the number of  buds of $T_e^\circ$ that are matched  
 to a  leaf of $T_e^\bullet$.   Let $j_{2,2}$ be
 the number of  buds of $T_e^\circ$ that are matched  to a leaf of
 $T_e^\circ$ in such a way the matching edge goes around the root of
 $T$ before ending at a leaf of $T_e^\circ$. Finally, let
 $j_2=j_{1,2}+j_{2,2}$. 
 In the closure of
 % {\st{$T_e^\circ$}
 $T_2$, the (black) root vertex  lies at
 depth  $j_2$ (Figure~\ref{fig:appendix},  bottom right).
Observe that, if the second single leaf of $T$
 is in $T_e^\bullet$, then $j_{2,2}=0.$
Since $T_e^\circ$ has charge $0$, counting its leaves and buds gives
$j_2= j_1-1 +j_{2,2}+s$, where $s\in\{1,2\}$ is the number of single
leaves of $T$ that lie in $T_e^\circ$. Consequently, $j_1 \le j_2$,
and %the
equality holds if and only if  the second single leaf of $T$,
denoted $\ell'$, is in $T_e^\bullet$.

% Finally define $j_1$ (resp.  $j_2$) to be the depth of
%the root of $T_1$ (resp.  $T_2$) in the matching sequence of its
% closure, considered as a bracket sequence.
%
% In the closure of $T$, no edge can wind around $T_e^\circ$: this
% subtree contains $\ell$, a single leaf of $T$.  Hence $j_1\leq j_2$.
% More over the second single leaf $\ell'$ of $T$ is in $T_e^\bullet$ if
% and only if $j_1=j_2$. Otherwise it is in $T_e^\circ$.  

If $\ell',\ell $ and $e$ appear in this order in counterclockwise
direction around   $T_e^\circ$, 
 set $ \rho(T)=(T_2,T_1)$.  In all other cases (and in particular if
$\ell'$ is in $T_e^\bullet$), set $ \rho(T)=(T_1,T_2)$.
Either ways, $ \rho(T)$ belongs to $(\mathcal{B}_2)^2$.
 Observe that if we write $\rho(T)=(T_a, T_b)$, then $a=2$ if and only if  $j_a>j_b$, with obvious notation.

\begin{figure}[htb]
\begin{center}
\epsfig{file=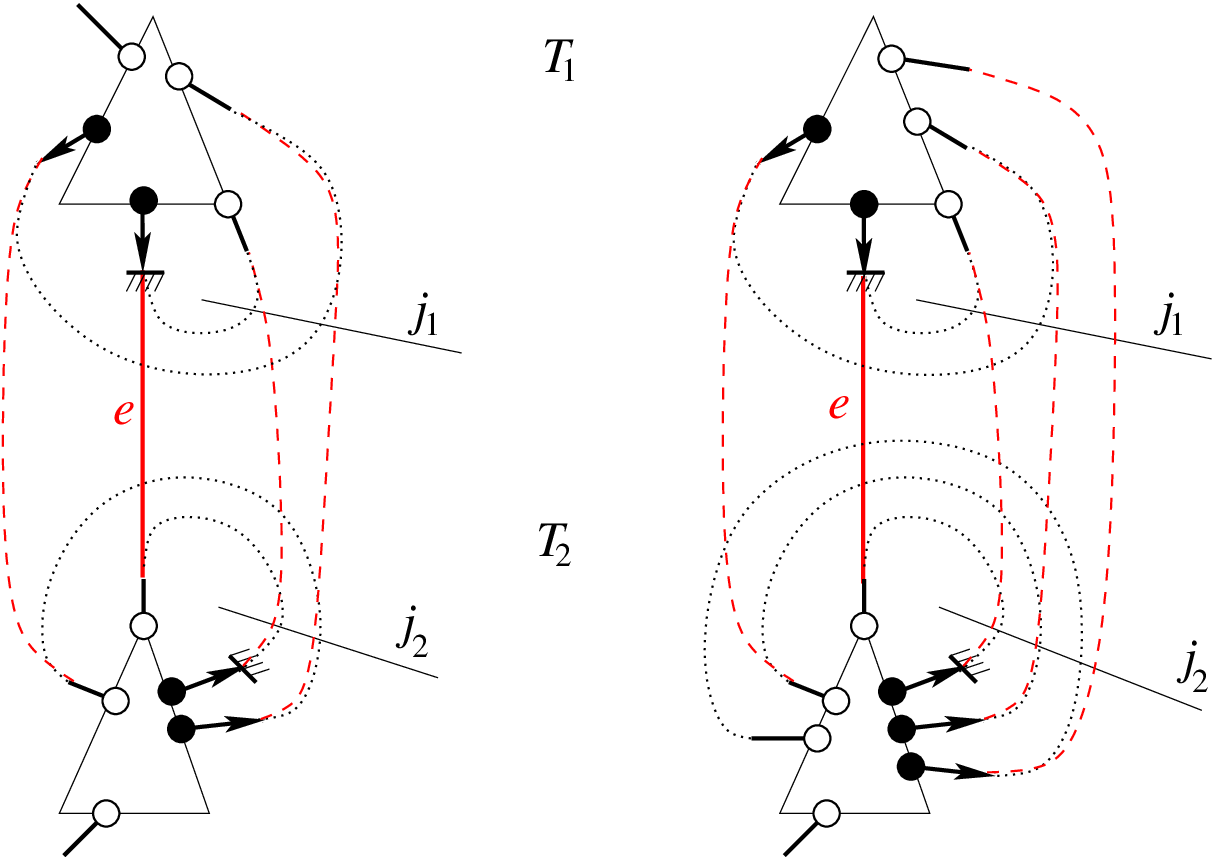, width=11cm}
\end{center}
\caption{Gluing two trees $T_1$ and $T_2$  of $\mathcal B_2$ such that $j_1\le j_2$. The dotted lines show the matchings of $T_1$ and $T_2$ (at least, the edges of the closure that are involved in the depths of the roots), while the (red) dashed lines show what they become after gluing the trees. The other matching edges do not change. On the left, $j_1=j_2$, while $j_1<j_2$ on the right.}  
\label{fig:leaves}
\end{figure}

Conversely, take a pair $(T_a,T_b)$ of $(\mathcal{B}_2)^2$. We will construct a tree from $T_a$ and $T_b$, which is \emm the only possible tree $T$ that is likely to be in $\mathcal W_\ell$ and to satisfy, $\rho(T)=(T_a, T_b)$. We will then observe that these two properties hold indeed. In other words, we claim that the reverse construction of~$\rho$ is forced.

First, note that $T_a$ and $T_b$ have total charge $1$. Let $j_a$
and $j_b$ be the depths of the roots of $T_a$ and~$T_b$ in the associated closures.
Rename the pair $(T_1,T_2)$ if $j_a\leq j_b$, and $(T_2,T_1)$ if $j_a>j_b$.
%otherwise.
The respective depths $j_1$ and $j_2$ of the roots of $T_1$ and $T_2$ then satisfy 
$j_1\leq j_2$.
Merge the root bud of $T_1$ with the leaf of $T_2$ matched (in $T_2$) to the root bud of
$T_2$. This creates an edge $e$ and
% an unrooted
a tree~$T'$ of total charge $2$, say unrooted for the moment. The white charge at~$e$ is $0$, and the black charge at $e$ is $2$. In the blossom trees $T_1$ and $T_2$, the white charge rule is satisfied at any edge by Lemma~\ref{lem:weakedges}. Now for an edge of $T'$ distinct from $e$, the white charge is either unchanged or increased by $1$ when constructing $T'$. This implies that  $T'$ satisfies the white charge rule at every
edge. By Lemma~\ref{lem:white-charge}, we can choose a root leaf for $T'$ that yields a blossom tree of $\mathcal W_1$, still denoted $T'$.

Let $C$ be the core of $T'$. Since $e$ is weak, the root of $T'$ is in $T^\bullet_e$, that is, in $T_1$.  Moreover, the black endpoint of $e$, say $v_\bullet$,  is in the core: indeed, if there was a weak edge of $T'$ on the path between the root of $T'$ and $v_\bullet$, then the tree $T_1$, whose root vertex is $v_\bullet$, would have a black-rooted lower subtree of charge $2$, and would not be a blossom tree.  Consider the
% two single leaves of
matching procedure in $T'$ (Figure~\ref{fig:leaves}). If $j_1=j_2$
then there is one single leaf of $T'$ in each of $T_1$ and $T_2$ (recall that they have total charge $1$);  see Figure~\ref{fig:leaves}, left. In this case let  $\ell$ be the single leaf of~$T_2$.
Otherwise,  $j_1< j_2$ and both
single leaves of $T'$  lie in $T_2$;  see Figure~\ref{fig:leaves}, right. If
the original pair $(T_a,T_b)$ is $(T_2,T_1)$, that is, if $j_a>j_b$, name these leaves $\ell$ and $\ell'$ so that the triple
$(\ell',\ell,e)$ appears in counterclockwise direction around $T_2$. Otherwise, name them so that $(\ell',\ell,e)$ appears in clockwise direction around $T_2$.

Define $\sigma$
and $\tau$ as above on the tree $T'$.
The single leaf $\ell$ is the endpoint of a $\sigma$-$\tau$ chain, and the %origin
starting point of this chain is a leaf
of the core.  Re-rooting $T'$ on this leaf yields a tree $T$ of $\mathcal W_1\setminus \mathcal W_1^*$. This tree belongs in fact to
$\mathcal{W}_\ell$, because  $\sigma$ and $\tau$ are 
independent of the choice  of the root (taken in the core).
Moreover, we can now check that 
$\rho(T)=(T_a,T_b)$.  This proves that $\rho$ is one-to-one
between $\mathcal{W}_\ell$ and $(\mathcal{B}_2)^2$, and concludes the
proof of Theorem~\ref{thm:bijection}.
\qed

%%%%%%%%%%%%%%%%%%%%%%%%%%%%%%%%%%%%%%%%%%%%%%%%%%%%%%%%
\section{Hard particle models on planar maps}
\label{sec:hard}
%%%%%%%%%%%%%%%%%%%%%%%%%%%%%%%%%%%%%%%%%%%%%%%%%%%%%%%%

In this section, we consider rooted planar maps in which some vertices
are occupied by a particle, in such a way two adjacent vertices are
not both occupied. In our figures, the vacant (resp.~occupied) vertices are
shown in white (resp.~black). Moreover, we require the two
ends of the root edge to be vacant.
We shall say %, for short,
that the map is \emm rooted at a vacant edge.,
Let $H(\bm X,\bm Y)$ denote the \gf \ of these maps, in
which $X_k$ (resp.~$Y_k$)  counts white (resp.~black) vertices of
degree $k$. As above, $\bm X$ stands for $(X_1,X_2, \ldots)$ and $\bm
Y$ for $(Y_1,Y_2, \ldots)$.

\begin{Theorem}\label{thm:hardmain}
  The hard particle \gf\ of planar maps rooted at a vacant edge 
can  be expressed in terms of the \gf \ of  blossom trees, defined by~\eqref{WBdef}--\eqref{equ:blackgf}, as follows: 
  \[
  H(\bm X,\bm Y) = (W_0-B_2)^2+W_1-B_3-B_2^2 ,
  \]
  where the series $W_i\equiv W_i(\bm x,\bm y)$ and $B_i\equiv B_i(\bm
  x,\bm y)$ 
  are evaluated at $x_k=X_k$ for $k\ge 1$, $y_2=1+Y_ 2$ and $y_k=Y_k$
  for $k \not 
  =2$.
\end{Theorem}

\begin{figure}[htb]
\begin{center}
\epsfig{file=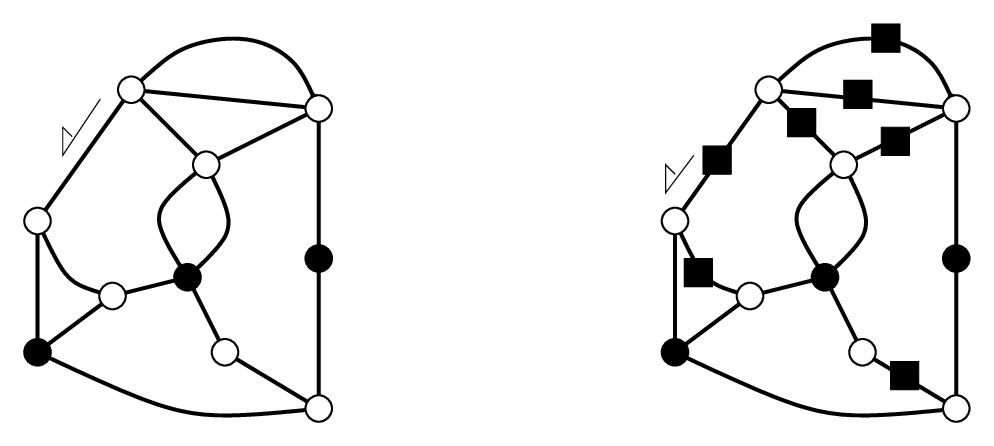, width=10cm}
\end{center}
\caption{From a hard particle configuration to  a bipartite map.}
\label{fig:hard-particle}
\end{figure}

\begin{proof} Take a map with hard particles. On every edge having both
ends vacant, add a black vertex of degree 2, of a special shape:
say, a square black vertex (Figure~\ref{fig:hard-particle}). One thus
obtains a bipartite map satisfying the following conditions:

-- the black vertices of degree 2 can be discs or squares,

-- all other vertices are discs,

-- the root vertex is a black square of degree 2.

\noindent We conclude using Theorem~\ref{thm:main}.
\end{proof}

% ============================================================
\subsection{Hard particles on tetravalent maps}
\label{section:hard4V}
%============================================================
We apply here Theorem~\ref{thm:hardmain} to tetravalent maps. That is,
all the variables $X_k$ and $Y_k$ are zero, except for $k=4$.

As a preliminary result, we need to
enumerate blossom trees having white vertices of degree~4 and black
vertices of degrees 2 and 4. We are actually going to solve a slightly
more general enumeration problem, by counting blossom trees having black
\emph{and} white vertices of degrees $2$ and $4$: first because this
 problem is nicely symmetric, then (and most importantly)
because we shall need this result later to solve the Ising model on
tetravalent maps. 

The corresponding series $W_i$ and $B_i$
 depend on the four variables $x_2$, $y_2$, $x_4$ and $y_4$, which
 we denote below by $v$, $w$, $x$ and $y$ to lighten notation.
 % for the sake of simplicity. Observe
 By Lemma~\ref{lem:multipleofm},
 the charge at the root of such blossom trees is always odd: hence the
series $W_{2i}$ and $B_{2i}$ are all zero. In particular, Theorem~\ref{thm:hardmain} reduces to $H(\bm X,\bm Y)=W_1-B_3$.
Equations~\eqref{equ:whitetrees} and~\eqref{equ:blacktrees} give
%specialize to 
\[
\left\{
\begin{array} {lll}
 W_1&=&v(1+B_1)+3xB_{-1}(1+B_1)^2, \\
W_3 &= &x(1+B_1)^3 ,
\end{array} 
\right.
\hskip 10mm 
\left\{
\begin{array} {lll}
  {B_{-3}} &=& {y}, \\
  B_{-1}&=& w+ 3yW_1, \\
B_1&=& wW_1+3y(W_3+W_1^2),
\end{array}
\right.
\]
 while $W_i=0$ for $i\ge 5$ and $B_i=0$ for $i\le -5$.
Moreover,~\eqref{equ:blackgf}  gives:
\beq \label{tree24-1}
B_3= wW_3+yW_1(6W_3+W_1^2).
\eeq
After a few reductions, we can express all the series $W_i$ and $B_i$
in terms of the series $P=1+B_1$, which satisfies: 
\beq
\label{tree24-2}
P=1+ 3xyP^3+ \frac{P(v+3xwP)(w+3yvP)}{(1-9xyP^2)^2}.
\eeq
In particular,
\beq \label{tree24-3}
W_1= \frac{P(v+3xwP)}{1-9xyP^2} \quad \hbox{and}\quad 
W_3= xP^3.
\eeq
By Theorem~\ref{thm:hardmain}, the hard particle \gf \ on tetravalent
maps is $W_1-B_3$,
%can be expressed in terms of the above series $W_1$ and $B_3$
with $v=0$ and $w=1$. We thus obtain the following result.

\begin{propo}[Hard particles on tetravalent maps]
The hard particle \gf \ of  tetravalent maps rooted at a vacant edge
is algebraic of degree $7$ and can be expressed as: 
\[
H(x,y)=xP^3+\frac{xP^2(3-2P)}{1-9xyP^2} -\frac{27x^3yP^6}{(1-9xyP^2)^3}
\]
where $P\equiv P(x,y)$ is the power series defined by
\[
P=1+ 3xyP^3 + \frac {3xP^2}{(1-9xyP^2)^2}.
\]
\end{propo}

\noindent{\bf Remarks}\\
{\bf 1.}  This proposition is exactly Corollary~\ref{cor:hard4V}.
In Section~\ref{subsec:Ising-hard} we had derived it from a  (more general) result dealing with the Ising model.\\
{\bf 2.} The parametrization by $P$ is equivalent the one given in
\cite{BDFG02a} for the free energy of this hard particle model: more
precisely, upon setting $x=g, y=zg$ and $P=V/g$, the equation defining
the parameter $P$ becomes Eq.~(2.14) of the above reference.

% ============================================================
\subsection{Hard particles on trivalent maps}
\label{section:hard3V}
%============================================================
We apply here Theorem~\ref{thm:hardmain} to trivalent (or cubic)
maps. That is, all the variables $X_k$ and $Y_k$
are zero, except for $k=3$. 

We  first need to
enumerate blossom trees having white vertices of degree 3 and black
vertices of degree 2 and 3. Again, we shall solve a more symmetric
problem by counting blossom trees having black \emph{and} white
vertices of degree $2$ and $3$. The corresponding series $W_i$ and $B_i$
 depend on  $x_2$, $y_2$, $x_3$  and $y_3$, which are denoted below by
$v, w, x$ and 
$y$.  
This model is a bit more complex than the previous one, since now
a blossom tree may have an even charge. However,
Equations~\eqref{equ:whitetrees} and~\eqref{equ:blacktrees} give
%specialize to 
\[
\left\{
\begin{array}  {lll}
W_0&=& vB_0+ x(2B_{-1}+2B_{-1}B_1+B_0^2), \\
W_1&=& v(1+B_1)+2xB_0(1+B_1), \\
W_2&=& x(1+B_1)^2,
\end{array}
 \right.
\hskip 5mm
\left\{
\begin{array}  {lll}
B_{-2}&=&y,\\
B_{-1}&=&w+2y  W_0, \\
B_0&=&wW_0+y(2W_1+W_0^2), \\
B_1&=&wW_1+2y(W_2+W_0W_1),
\end{array}
 \right.
\]
while $W_i=0$ for $i\ge 3$ and $B_i=0$ for $i\le -3$.
Equation~\eqref{equ:blackgf} then gives
\beq
\label{tree23-1}
B_2=wW_2+y(2W_0W_2+W_1^2)\quad \hbox{and} \quad B_3=2yW_1W_2 .
\eeq
Let  $P=1+B_1$, $Q=W_0$ and $R=B_0$. Then
\beq
\left\{ \begin{array}{l}
\displaystyle P =  1 + 2xy P^2 + P(w+2yQ)(v+2xR) , \\
\\
\displaystyle (1-4xyP) Q = 2xwP + R(v+xR) ,\\
\\
\displaystyle (1-4xyP) R =  2yvP + Q(w+yQ) ,
\end{array}\right.
\label{tree23-2}  
\eeq
and all the series $W_i$ and $B_i$ have simple rational expressions in
terms of $P,Q,R$. 
In particular,
\beq \label{tree23-3}
W_0=Q, \quad  B_{-1}=w+2yQ,\quad W_1= P(v+2xR) \quad \hbox{and}\quad W_2=xP^2.
\eeq 
By Theorem~\ref{thm:hardmain}, the hard particle \gf \ on trivalent maps
can be expressed in terms of the above series $W_i$ and $B_i$
evaluated at $v=0$ and $w=1$. We thus obtain the following result.

\begin{propo}[Hard particles on trivalent maps]
  \label{prop:hard3}
  The hard particle \gf \ of 
trivalent maps rooted at a vacant edge  is algebraic of degree $11$
and can be expressed as follows: 
\[
H(x,y)=Q^2+2xPR-4x^2yP^3R-2xP^2Q-4xyP^2Q^2-8x^2yP^2QR^2,
\]
where $P,Q$ and $R$  are the power series defined by
Eqs.~{\em\eqref{tree23-2}} above, with $v=0$ and $w=1$.
\end{propo}

\begin{proof}
 The expression of $H$ comes directly from Theorem~\ref{thm:hardmain} and the equations listed above the proposition. For the degree, we proceed as follows.
  Set $v=0$ and $w=1$ in~\eqref{tree23-2}, and eliminate~$R$ using the third equation. This gives
\beq\label{param-HP3}
\left\{
\begin{array}{lll}
P&=&\displaystyle 1+2xyP^2+\frac{2xPQ(1+yQ)(1+2yQ)}{1-4xyP},\\
\\
Q&=&\displaystyle  \frac{2xP}{1-4xyP}+ \frac{xQ^2(1+yQ)^2}{(1-4xy  P)^3}.
\end{array}
\right.
\eeq
  We then eliminate $R$, $P$ and $Q$ from the expression of $H$ to obtain the minimal polynomial of~$H$. The series $P$, $Q$, and $R$ also have degree $11$ and lie in the same extension of degree $11$ of $\qs(x,y)$ as $H$.
\end{proof}

\noindent{\bf Remark.} The algebraic parametrization~\eqref{param-HP3} is easily ckecked to be the same as the
one given in~\cite{BDFG02a} for the free energy of the
hard particle model on trivalent maps: more precisely, upon setting
$x=g, y=gz, P=V/(zg^2)$ and $Q=R/(zg)$, one recovers Eqs.~(C.7)
and~(C.8) of the above reference. 

%%%%%%%%%%%%%%%%%%%%%%%%%%%%%%%%%%%%%%%%%%%%%%%%%%%%%%%%
\section{The Ising model on planar maps}
\label{sec:ising}
%%%%%%%%%%%%%%%%%%%%%%%%%%%%%%%%%%%%%%%%%%%%%%%%%%%%%%%%
In this section, we consider maps with white and black vertices, and
enumerate them according to their degree distribution (variables $X_k$
and $Y_k$) and according to the number of frustrated edges (variable
$u$). Let us recall that an edge is frustrated if it has a black endpoint
and a white one. The resulting series is called the \emm Ising \gf, of the class of maps under consideration.

% ============================================================
\subsection{General result}
\label{section:Ising-general}
% ============================================================

\begin{Theorem}\label{thm:Ising-general}
The Ising \gf \ $I(\bm X, \bm Y, u)$ of planar maps  rooted at a black vertex
of degree $2$ can be expressed in terms of the \gfs \ for blossom trees, defined by~\eqref{WBdef}--\eqref{equ:blackgf}, as follows: 
\[
I(\bm X, \bm Y, u)= Y_2(u-\bu) \left( (W_0-B_2)^2+W_1-B_3-B_2^2 \right),
\]
where the series $W_i\equiv W_i(\bm x,\bm y)$ and $B_i\equiv B_i(\bm x,\bm y)$
  are evaluated at 
\[\left\{\begin{array}{lllllll}
x_2&=&\bu + X_2(u-\bu), \\
x_k&=&X_k(u-\bu)^{k/2} & \hbox{for } k \not =2,
         \end{array}\right.
       \hskip 8mm
\left\{\begin{array}{lllllll}
 y_2&=&\bu + Y_2(u-\bu), \\
 y_k&=&Y_k(u-\bu)^{k/2}& \hbox{for } k \not =2,
\end{array}\right.
\]
with $\bu =1/u$.
\end{Theorem}

\begin{proof} Take a bicoloured map rooted at a black vertex of
degree~2. On each edge, add a (possibly empty) sequence of square
vertices of degree 2, in such a way the resulting map is bipartite. An
example is shown in Figure~\ref{fig:ising-bip}. Note that every
frustrated edge receives an even number of square vertices, while
every monochromatic edge receives an odd number of these squares. The
resulting map remains rooted at its black vertex of degree 2, which is
not a square.

\begin{figure}[htb]
\begin{center}
  \epsfig{file=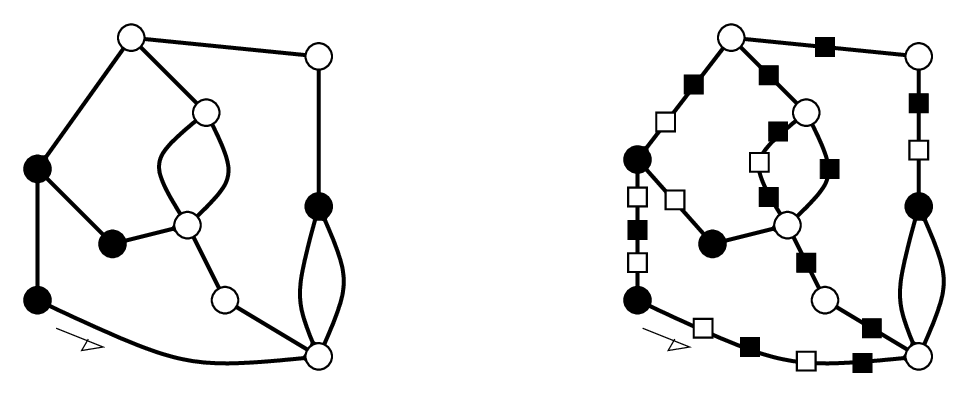, width=10cm}
\end{center}
\caption{An Ising configuration with a black root of degree two and one of the associated bipartite maps.}
\label{fig:ising-bip}
\end{figure}

Let $\widetilde M(\bm x,\bm y,v)$ denote the degree \gf \ of  the maps one
obtains in that way: in this series, the square vertices are counted
by $v$, while the other vertices are, as usually, counted by the
variables $x_k$ and $y_k$. Given that the total number of edges in the original map is $\frac 1 2 \sum_k k(m_k+n_k)$, where $m_k$ (resp.~$n_k$) is the number of white (resp. black) vertices of degree~$k$, the above construction gives
\beq
\widetilde M(\bm x,\bm y,v)= I(\bm X,\bm Y,u)
\label{Mtilde1}
\eeq
where 
\[
X_k = \frac{x_k}{(\bv -v )^{k/2}}, \quad Y_k = \frac{y_k}{(\bv -v)^{k/2}}
\quad \hbox{and} \quad u=1/v =\bv.
\]
By Theorem~\ref{thm:main}, we also have
\beq
\widetilde M(\bm x,\bm y ,v)= y_2 \left( (W_0-B_2)^2+W_1 -B_3 -B_2^2\right)
\label{Mtilde2}
\eeq
where the series $W_i$ and $B_i$ are evaluated at $x_1, x_2+v, x_3,
\ldots, y_1, y_2+v, y_3, \ldots$ The result follows by comparison
of~\eqref{Mtilde1} and~\eqref{Mtilde2}.
\end{proof}

%============================================================
\subsection{The Ising model on quasi-regular maps}
%============================================================
The general result above applies in particular to quasi-regular maps,
that is, maps
in which all vertices have degree $m\ge 3$, except for the root which has
degree $2$. Let us make Theorem~\ref{thm:Ising-general} explicit for
the case $m=4$ 
(this is Proposition~\ref{propo:Isingquasi4V} in Section~\ref{sec:glimpse}) and $m=3$
(Proposition~\ref{thm:Isingquasi3V} below). 

\begin{proof}[Proof of Proposition~\ref{propo:Isingquasi4V}]
  % For  quasi-tetravalent maps,
  We first form the Ising
\gf \ of  planar maps, rooted at a black vertex of degree $2$, having
only vertices of degree $2$ and $4$. This series is given by
Theorem~\ref{thm:Ising-general}, with 
$X_k=Y_k=0$ for $k \not = 2,4$.
The blossom trees
% occurring in
underlying  this theorem are the ones that have been counted in Section~\ref{section:hard4V}. Their enumeration has resulted in  
Eqs.~\eqref{tree24-1}--\eqref{tree24-3}.
Theorem~\ref{thm:Ising-general} yields:
\[
I(\bm X, \bm Y, u)= Y_2(u-\bu) (W_1-B_3),
\]
where $W_1$ and $W_3$ are evaluated at $v= \bu + X_2(u-\bu)$,  $w=
 \bu + Y_2(u-\bu)$, $x=X_4(u-\bu)^2$, $y=Y_4(u-\bu)^2$.

 The Ising \gf \ $\quasiI(X,Y,u)$ defined in Proposition~\ref{propo:Isingquasi4V}
 % of   quasi-tetravalent maps
 is obtained by setting
 $X_2=0 $ in the series $I(\bm X, \bm Y, u)$, then extracting
 the coefficient of $Y_2$, and finally replacing~$X_4$ and $Y_4$ by $X$ and $Y$, respectively. This gives $\quasiI(X,Y,u)=(u-\bu) (W_1-B_3)$ where $W_1$ and~$W_3$ are evaluated at $v=w=\bu$, $x=X(\bu-u)^2$, $y=Y(\bu-u)^2$.
Proposition~\ref{propo:Isingquasi4V} follows using~\eqref{tree24-1} and~\eqref{tree24-3}.
\end{proof}

We now address the case of quasi-cubic maps.

\begin{propo}[Ising on quasi-cubic maps]
  \label{thm:Isingquasi3V}
Let $\quasiI( X,Y,u)$ be the Ising \gf \ of quasi-cubic maps rooted
at a black vertex, where
the variables $X$
and $Y$  account for the number of (cubic) white and 
black vertices, while $u$ counts  frustrated edges.
Then
\[
\frac{\quasiI(X,Y,u)}{u-\bu}=Q^2+P(v+2xR)-2xyP^3(v+2xR)
 -2xvP^2Q
-4xyP^2Q^2-2yP^2Q(v+2xR)^2,
\]
where  $P$, $Q$, $R$ are the series  defined by~\eqref{tree23-2},
evaluated at 
 $x=X(u -\bu)^{3/2}$,  $y=Y(u -\bu)^{3/2}$ 
and $v=w=1/u=\bu$. All these series are algebraic of degree $11$. 
\end{propo}

\begin{proof}%[Proof of Proposition~\ref{thm:Isingquasi3V}]
  The  argument used above in the proof of Proposition~\ref{propo:Isingquasi4V}
  gives
  %the Ising \gf \ of quasi-cubic maps as
\[
\quasiI(X,Y,u)= (u-\bu)\left((W_0-B_2)^2+W_1-B_3-B_2^2\right),
\]
where the series $W_i$ and $B_i$ are given by
Eqs.~\eqref{tree23-1}--\eqref{tree23-3} and are evaluated at $v=w=\bu$,
$x=X(u-\bu)^{3/2}$, $y=Y(u-\bu)^{3/2}$. The degree of $\quasiI$ is established by elimination of $P$, $Q$ and $R$, as in the proof of Proposition~\ref{prop:hard3}.
\end{proof}

\noindent{\bf Remark.}
%{\bf 1.}
Set $w=v$ in~\eqref{tree23-2}, and denote $\bar Q=
v+2yQ$, $\bar R = v+2xR$. Then
\[
\left\{
\begin{array}{l}
\displaystyle P =  1 + 2xy P^2 + P\bar Q\bar R , \\
\\
\displaystyle (1-4xyP) \bar Q = v + \frac y{2x} \left(\bar R^2-v^2\right) ,\\
\\
\displaystyle (1-4xyP) \bar R = v + \frac x{2y} \left(\bar Q^2-v^2\right) .
\end{array}
\right.
\]
The above system of equations is equivalent to the
parametrization given in~\cite{BK87} for the free energy of the Ising
model on trivalent maps. More precisely, upon setting
\[
x= \frac{ g e^H} {c^{3/2}},\quad y= \frac{ g e^{-H}} {c^{3/2}}, \quad v=1/c,
\quad
P=\frac{c^2 z}{g^2}, 
\quad \bar Q= -\frac {\kappa}{c}, \quad  \bar R = -\frac {\rho}{c},
\]
the parametrization in $P,\bar Q,\bar R$ becomes Eqs.~(40--41) of
Ref.~\cite{BK87}, with $g^2$ replaced by $-g^2$.\\
%

% ============================================================
\subsection{The Ising model on regular maps}
% : integral solution}
\label{sec:Ising-regular}
%============================================================
Erasing the root vertex of a quasi-regular map gives an authentic
regular map. This implies that
% anks to this remark, the above
the \gfs \ denoted $\quasiI$ in Propositions~\ref{propo:Isingquasi4V} and~\ref{thm:Isingquasi3V} can also be
written as
\[
\quasiI(X,Y,u)= u^2 I_{0,0}(X,Y,u) + I_{0,1}(X,Y,u) + I_{1,0}(X,Y,u) +
I_{1,1}(X,Y,u) , 
\]
where $I_{0,0}$ is the Ising \gf \ of  $m$-regular maps rooted at a
% uniformly
white edge, $I_{0,1}$ is the Ising \gf \ of  $m$-regular
maps rooted at a frustrated edge (oriented from its white endpoint to
its black one) and so on (in Propositions~\ref{propo:Isingquasi4V} and~\ref{thm:Isingquasi3V}, $m$ is $4$ and $3$ respectively).
However, this does not give us
%The general result of Theorem~\ref{thm:Ising-general} is obviously not
%very convenient for solving
the Ising model on truly regular
$m$-valent maps. We will show  that a re-rooting procedure circumvents this difficulty,
% to the cost of an
but results in an integral expression.
% expression of the Ising \gf \ in terms of an integral.
In the tetravalent case at least, we will evaluate this integral  as an
algebraic {series.} %function.

\begin{figure}[htb]
\begin{center}
  \epsfig{file=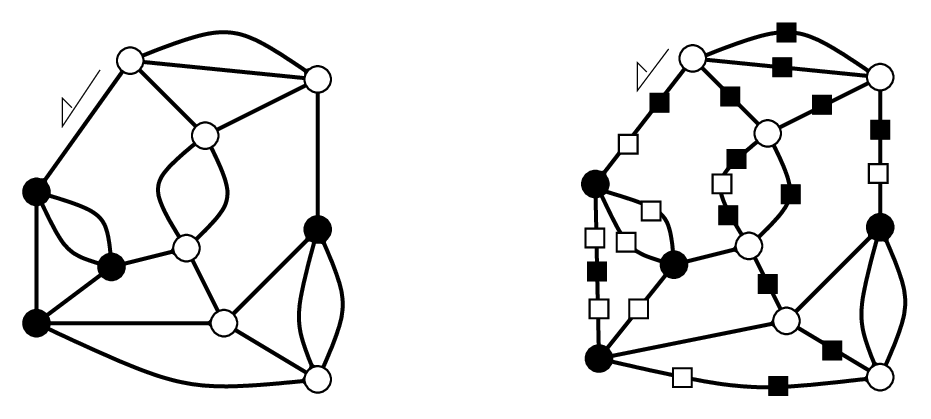, width=10cm}
\end{center}
\caption{An Ising configuration on a tetravalent map and one of the associated bipartite maps.}
\label{fig:ising-tetra}
\end{figure}

Suppose we are interested in the Ising enumeration of $m$-valent maps,
rooted at a white vertex, with $m\ge 3$.  We still denote by
$I(X,Y,u)$ the associated
\gf , where $X$ (resp.~$Y$) counts white (resp.~black) vertices, and
$u$ counts  frustrated edges.  As in Section~\ref{section:Ising-general}, let us add a sequence of vertices of degree $2$ on every edge, in such a way
the resulting map is bipartite, and still rooted at a white vertex of degree $m$ (Figure~\ref{fig:ising-tetra}). Let
$\widetilde M (x,y,v)$ be the \gf \ of  the bipartite maps thus obtained, where
the  vertices of degree 2 are counted by $v$, and the other ones by
$x$ and $y$ depending on their colour. Then
\beq\label{Mising0}
\widetilde M (x,y,v)= I \left( \frac x{(\bv -v )^{m/2}}, \frac y{(\bv -v
)^{m/2}}, \bv\right).
\eeq
We now explain how to relate $\widetilde M$ to the series $W_i$ and $B_i$.

First,
% Now,
let $\overline M(x,y,v)$ be the \gf \
% $M(\bm x,\bm y)$
of bipartite maps
rooted at a black vertex of degree 2 (given by Theorem~\ref{thm:main}),
evaluated at $x_2= y_2=v, x_m=x, y_m=y$, all the other variables
 $x_k$ and $y_k$ 
being zero. Then $\overline M(x,y,v)+\overline M(y,x,v)$ counts
maps rooted at a vertex of degree 2, either black or white.
In these maps, let us orient another edge,
starting now from a white vertex of degree $m$: we get a doubly rooted
map, which can be obtained in an alternative may, namely by orienting
an edge that starts from a vertex of degree 2 in a map rooted at a white vertex of degree $m$ (these are the maps counted
by $\widetilde M (x,y,v)$). In algebraic terms, this yields
\[
mx \frac{\partial}{\partial x} \left( \overline M(x,y,v)+\overline M(y,x,v)\right) =
 2v \frac{\partial
\widetilde M}{\partial v} (x,y,v)  .
\]
Integrating this equation with respect to $v$ provides the \gf \ of bipartite maps rooted
at a white vertex of degree $m$ that contain at least one vertex of
degree $2$:
\beq
\widetilde M(x,y,v)-\widetilde M(x,y,0)= \frac {mx} 2 \int_0^v
\frac{\partial}{\partial x} \left( \overline M(x,y,z)+\overline M(y,x,z)\right)
 \frac {dz} z.
\label{Mising1}
\eeq
Since $\overline M$ is given by Theorem~\ref{thm:main}, what remains is to evaluate the \gf \ $\widetilde M(x,y,0)$ of $m$-valent bipartite
maps rooted at a white vertex. By deleting the black endpoint of the root 
edge, together with the $m$ half-edges that are attached to it, we
obtain an $m$-leg map, as defined in Section~\ref{sec:legs}. By Theorem~\ref{thm:closure-m}, these maps are
in bijection 
with  $m$-regular balanced blossom trees of total charge $m$. For
$m$-regular blossom trees, Eqs.~\eqref{equ:whitetrees}
and~\eqref{equ:blacktrees} simply give: 
\beq \label{equ:m-valent}
W_{m-1}=x(1+B_1)^{m-1} \quad \hbox{and} \quad B_1=(m-1)yW_{m-1},
\eeq
while Eq.~\eqref{equ:blackgf} gives 
\[
B_{m+1}=y {{m-1} \choose 2} W_{m-1}^2.
\]
Now Theorem~\ref{thm:constellationplusplusbifluore} gives the \gf \
of $m$-regular balanced blossom trees of total charge $m$ as $W_{m-1}-B_{m+1}$.  By Theorem~\ref{thm:closure-m},  the \gf\ of $m$-valent bipartite
maps rooted at a white vertex is then
\beq
\widetilde M(x,y,0)=  y \left(W_{m-1}-B_{m+1}\right) =A - {{m-1}\choose 2} A^2,
\label{Mising3}
\eeq
where $A=yW_{m-1}$ is the only formal power series satisfying
\beq
\label{equ:constellations}
A=xy(1+(m-1)A)^{m-1}.
\eeq
By combining~\eqref{Mising0},~\eqref{Mising1} and~\eqref{Mising3}, we obtain
an integral expression of the Ising \gf \ of $m$-valent maps.

\begin{Theorem}\label{thm:Ising-regulier}
Let $m\ge 3$. The Ising \gf \ $I(X,Y,u)$ of $m$-valent maps rooted at a white vertex can be expressed as follows:
\[
I(X,Y,u)=  A - {{m-1}\choose 2} A^2
+ \frac {mx} 2 \int_0^v
\frac{\partial}{\partial x} \left( \overline M(x,y,z)+\overline M(y,x,z)\right)
 \frac {dz} z,
\] 
where the series $A\equiv A(x,y)$ is given by~{\em
\eqref{equ:constellations}} and $\overline M(x,y,z)$ counts bipartite maps
having only vertices of degree $2$ and $m$, rooted at a black vertex
of degree~$2$. This series has an algebraic expression given by Theorem~\ref{thm:main}. Both series $A$ and $\overline M$ are evaluated at $x=X(u-\bu)^{m/2}$,
$y=Y(u-\bu)^{m/2}$ and $v=1/u=\bu$.
\end{Theorem}

% ============================================================
%\subsection{The Ising model on regular 4-valent maps: algebraicity}
%============================================================
As an application of the above general result, we shall now solve the Ising model
on truly tetravalent maps.
We shall see that, in this case at least,  the series $I(X,Y,u)$
 is itself algebraic, and belongs to the same extension of
 $\qs[X,Y,u]$ as
$\overline M(x,y,v)$. We expect this to be  always
 true: the nature of a map \gf \ should not depend on details of the
 choice of the root.

 \begin{propo}[Ising on tetravalent maps]\label{propo:Ising4V}
   Let $I(X,Y,u)$ be the
Ising \gf \ of  bicoloured tetravalent maps rooted at a white vertex.
Let $P(x,y,v)$ be the {algebraic} power series defined by the following
%algebraic
 equation:
  \beq\label{P-alg}
  P=1+3xyP^3 + v^2 \frac{P(1+3xP)(1+3yP)}{(1-9xyP^2)^2}.
  \eeq
  Then 
$I(X,Y,u)$ can be expressed in terms of
  $P(x,y,v)$, with $x=X(u- \bu)^2, y=Y(u-\bu )^2$ and $v= \bu =
  1/u$. One possible expression is:
  % \[I(X,Y,u)=  \frac 1 9 \Big(135x^2y^2P^6+72xy^2P^5 -3y(15x+8y-3v^2y+36xy)P^4\hskip 50mm\]\[\hskip 30mm  -y(32+3v^2-36x -36y)P^3 -(3+2v^2-72y)P^2+12(1-3y)P -9 \Big)\]\[\hskip 80mm +\frac{(1-P-yP^2)(1+3yP)(12-8P+ 3v^2P)}{9(1+3xP)}.\]
  \[
  I(X,Y,u)= xP \, \frac{\Pol(x,y,P)}{(1+3xP)^2},
\]
with
\begin{multline*}
\Pol(x,y,P)=  -243 x^{3} y^{4} P^{9}+162 x^{3} y^{3} P^{8}+135 x^{2} y^{2} \left(x +2 y \right) P^{7}-18 x^{2} y^{2} \left(9 y -4\right) P^{6}\\
  -3 x y \left(18 x y +15 x +19 y \right) P^{5}+12 x y \left(3 x +6 y -4\right) P^{4}-3 \left(9 x \,y^{2}-28 x y +x +y \right) P^{3}\\
  -2 \left(18 x y -6 x -5 y \right) P^{2}-\left(9 x +6 y -2\right) P -2.
\end{multline*}
As $P$ itself, the series $I(X,Y,u)$ is algebraic\footnote{But the
algebraic equation satisfied by $I$ is not  
something one likes to see.} of degree $7$.
\end{propo}

\begin{proof} 
We start from Theorem~\ref{thm:Ising-regulier}, taken with $m=4$. In order to express the series $\overline M$,  we have to count blossom trees having only vertices of
 degree $2$ and $4$. This has been done in
 Section~\ref{section:hard4V},  and has resulted in
 Eqs.~\eqref{tree24-1}--\eqref{tree24-3}.  We only need the case $w=v$,
 which gives the above parametrization for $P$. Note that it is symmetric in $x$ and $y$.

We shall now merely sketch the rest of the computation, for which it
is useful to use a computer algebra software, like {\sc Maple}. From
Theorem~\ref{thm:main} and Section~\ref{section:hard4V}, we obtain an expression of the \gf \ $\overline M(x,y,v)=v(W_1-B_3)$ that counts bipartite maps rooted at a black vertex of
degree~$2$ in terms of the above parameter $P$. By exchanging $x$ and~$y$, we obtain a similar expression for $\overline M(y,x,v)$. We
then compute the derivatives of $P$ with respect to~$x$ and to $v$. Both
derivatives are expressed as rational functions of $x,y,v$ and~$P$.  
Observe that, from the equation defining $P$, we can express
 any even function of $v$  in terms of $x,y$ and~$P$.

However, the integrand occurring in Theorem~\ref{thm:Ising-regulier}, now expressed rationally in terms of $x$, $y$, $v$ and $P$,
is an odd function of $v$ (or more precisely~$z$). But so is the $v$-derivative of $P$. The combination of these two remarks allows us to write
the integrand in the form 
\[ 
\frac{N(x,y,P)}{D(x,y,P)} \ \frac{\partial P}{\partial z}(x,y,z),
\]
where $N$ and $D$ are polynomials in $x,y$ and $P$, not involving
$z$. Hence the integral becomes 
\[
\int_{P(x,y,0)}^{P(x,y,v)}  \frac{N(x,y,p)}{D(x,y,p)} dp.
\]
The integrand is now an explicit rational function, and its
% primitive
antiderivative in $p$
turns out to be rational too. Hence the integral of
Theorem~\ref{thm:Ising-regulier} is a rational function of
$x,y, P(x,y,v)$ and $P(x,y,0)$. But the series $P(x,y,0)$ is closely
related to the series $A$: 
indeed, it follows from~\eqref{equ:constellations} and~\eqref{P-alg} that
  % $P(x,y,0)-1\equiv B_1(x,y,0)$ counts blossom trees with regular $4$-valent vertices, so that
  $P(x,y,0)=1+3A$ (this can be explained combinatorially as well).
  % (see~\eqref{equ:m-valent}).
Thanks to the equation~\eqref{equ:constellations} defining $A$,   the terms involving $A$ %miraculously
  vanish in the expression of $I(X,Y,u)$, and
we end up with the not-so-simple, but algebraic, formula of
Proposition~\ref{propo:Ising4V}. 
\end{proof}

\noindent{\bf Remarks} \\
{\bf 1.} The equation~\eqref{P-alg} defining $P(x,y,v)$ is equivalent to the
parametrization given in~\cite{BK87} for the free energy of the Ising
model on tetravalent maps. More precisely, upon setting
\[
x= \frac{ -g e^H} {c^{2}},\quad y= \frac{ -g e^{-H}} {c^{2}}, \quad v=1/c
\quad \hbox{and}\quad 
P=-{c^2 z}/(3g),\]
the above parametrization in $P$ becomes Eq.~(17) of
Ref.~\cite{BK87}.

\begin{figure}[htb]
\begin{center}
\epsfig{file=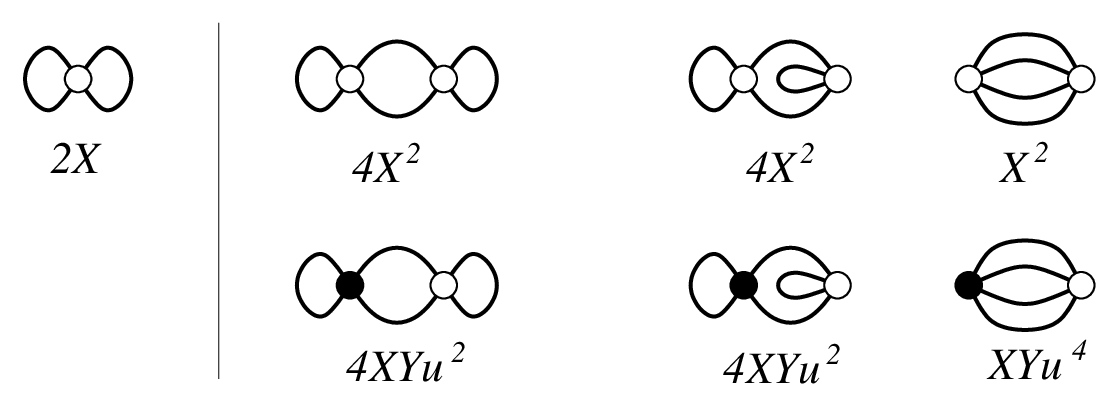, width=12cm}
\end{center}
\caption{The (unrooted) tetravalent maps with at least one white vertex
(and  a most two
vertices). The multiplicities give the number 
of ways to root at a white vertex, and each rooted map is then
weighted as in $I(X,Y,u)$.}
\label{fig:FirstOnes4V}
\end{figure}

\noindent
{\bf 2.} Set $X=tX$ and $Y=tY$, so that $t$ counts %the number of
vertices. Then
\[
I(tX,tY,u)= t(2X)+t^2(9X^2+XY(8u^2+u^4))+ \LandauO(t^3).
\]            
The corresponding Ising configurations are shown on
Figure~\ref{fig:FirstOnes4V}. 

%\small

\bibliographystyle{abbrv}
\bibliography{ising.bib}

\end{document}